\documentclass[10pt,twosided]{article}
\usepackage[tbtags]{amsmath}
\allowdisplaybreaks[4]
\usepackage{amssymb,color}

\definecolor{c20}{rgb}{0.,0.7,0.}
\definecolor{c30}{rgb}{0.,0.,1.}
\definecolor{c40}{rgb}{1,0.1,0.7}
\definecolor{c50}{rgb}{1,0,0}
\definecolor{c60}{rgb}{1,0.9,0.1}

\def\re#1{\textcolor{c30}{#1}}
\def\re#1{#1}

\def\cH#1{\textcolor{c20}{#1}}
\def\cH#1{#1}

\def\rL#1{\textcolor{c40}{#1}}
\def\rL#1{#1}

\newcommand{\kb}[1]{\boldsymbol{#1}}
\newcommand{\vk}[1]{\kb{#1}}

\newcommand{\ve}{\varepsilon}

\newcommand{\abs}[1]{\lvert #1 \rvert}

\newcommand{\E}[1]{\mathbb{E}\left\{#1\right\}}

\newcommand{\pk}[1]{\mathbb{P} \left\{ #1 \right\} }

\newcommand{\R}{\!I\!\!R}

\newcommand{\N}{\!I\!\!N}
\newcommand{\inr}{\in \R}

\newcommand{\inn}{\in \N}

\newcommand{\ldot}{,\ldots,}

\newcommand{\BQN}{\begin{eqnarray}}
\newcommand{\EQN}{\end{eqnarray}}
\newcommand{\BQNY}{\begin{eqnarray*}}
\newcommand{\EQNY}{\end{eqnarray*}}

\newcommand{\BS}{\begin{sat}}
\newcommand{\ES}{\end{sat}}
\newcommand{\BT}{\begin{theo}}
\newcommand{\ET}{\end{theo}}
\newcommand{\BK}{\begin{korr}}
\newcommand{\EK}{\end{korr}}

\newcommand{\BD}{\begin{de}}
\newcommand{\ED}{\end{de}}

\newcommand{\BIT}{\begin{itemize}}
\newcommand{\EIT}{\end{itemize}}
\newcommand{\BDI}{\begin{description}}
\newcommand{\EDI}{\end{description}}

\newcommand{\BRM}{\begin{remarks}}
\newcommand{\ERM}{\end{remarks}}

\newcommand{\BEL}{\begin{lem}}
\newcommand{\EEL}{\end{lem}}

\newcommand{\BLE}{\begin{lem}}
\newcommand{\ELE}{\end{lem}}

\newtheorem{theo}{Theorem}[section]
\newtheorem{sat}[theo]{Proposition}
\newtheorem{de}[theo]{Definition}
\newtheorem{lem}[theo]{Lemma}

\newtheorem{example}[theo]{Example}
\newtheorem{korr}[theo]{Corollary}

\newtheorem{remarks}[theo]{Remarks}

\newcommand{\nelem}[1]{{Lemma \ref{#1}}}
\newcommand{\neprop}[1]{{Proposition \ref{#1}}}
\newcommand{\netheo}[1]{{Theorem \ref{#1}}}

\newcommand{\prooftheo}[1]{ \textsc{Proof of Theorem} \ref{#1} }

\newcommand{\prooflem}[1]{\textsc{Proof of Lemma} \ref{#1}}

\newcommand{\COM}[1]{}

\newcommand{\QED}{\hfill $\Box$}

\def\rw{\rightarrow}

\date{}
\def\rw{\rightarrow}
\def\IF{\infty}
\def\LT{\left}
\def\RT{\right}
\def\det{\triangle}
\def\Ha{ \mathcal{H}_{\alpha}}
\def\P{ \mathcal{P}}
\def\vn{\varepsilon}
\def\On{\mathcal{S}_{n-1}}

\def\oo{o}
\def\ooo{(1+o(1))}
\def\vara{\varsigma}
\def\NN{\mathcal{N}_*}
\def\wtre{\widetilde{\triangle}}
\def\CK{\mathcal{K}}
\newcommand{\vkt}[1]{\kb{#1}}
\def\MM{\mathcal{M}}

\def\QQ{\mathcal{Q}}
\def\DD{\mathcal{D}}
\def\WW{\mathcal{W}}

\def\parO{\partial \mathcal{S}}
\def\pO{\partial \mathcal{O}}
\def\xih{\Bigl(\chi_n(t)- \rL{g(t)} \Bigr)}

\def\wsig{\widetilde{\sigma}}
\def\chicb{\rL{(\chi_n(t)-ct^\beta)}}

\begin{document}

\title{\bf \Large Piterbarg Theorems for Chi-processes with Trend}

\author{Enkelejd Hashorva and   Lanpeng Ji
\thanks{University of Lausanne, UNIL-Dorigny 1015 Lausanne, Switzerland}}

 \maketitle
  \centerline{\today{}}

{\bf Abstract:} Let  $\chi_n(t) = (\sum_{i=1}^n X_i^2(t))^{1/2},\  {t\ge0}$  be a
chi-process with $n$ degrees of freedom  where   $X_i$'s
are independent copies of some generic centered Gaussian process $\re{X}$. This paper derives the exact asymptotic behavior of
\BQNY
\pk{\sup_{t\in[0,T]} \xih >u}\ \ \ \text{as}\ \ u \rw \IF,
\EQNY
where \rL{$T$ is a given positive constant, and $g(\cdot)$ is some non-negative bounded measurable function.}
The case $g(t)\equiv0$ is investigated in numerous contributions by V.I. Piterbarg. Our novel asymptotic results,
for both \re{stationary and non-stationary $X$}, are referred to as Piterbarg theorems for chi-processes with trend.

{\bf Key words:}  Gaussian random fields; Piterbarg theorem for chi-process;  Pickands constant; generalized Piterbarg  constant;
 Piterbarg inequality.

{\bf AMS Classification:} Primary 60G15; Secondary 60G70.

\section{Introduction}

Two fundamental results for the study of asymptotic behaviour of supremum of \rL{non-smooth} Gaussian processes and Gaussian random fields are
{\it Pickands theorem} and {\it Piterbarg theorem}, see Pickands (1969a,b), Piterbarg (1972, 1996), and Piterbarg and Prisyazhnyuk (1978). 
\rL{For any fixed $T\in(0,\IF)$,}
J. Pickands III  obtained the exact asymptotics of the probability $\pk{\sup_{t\in[0,T]}X(t)>u}$ as $u\rw\IF$ for
a centered stationary Gaussian process $\{X(t),t\ge 0\}$
 with a.s. continuous sample paths and \rL{covariance} function \rL{$r(\cdot)$} satisfying the following  assumptions:

{\bf Assumption R1}.  $r(t) = 1 -\rL{\abs{t}^{\alpha}}\ooo$ as $t \rightarrow 0,$ with $\alpha \in (0, 2]$; 

{\bf Assumption R2}. $r(t)<1$ for all $t>0$.

More precisely, {\it Pickands theorem} states \re{that}
\BQN\label{eq:PP}
\pk{\sup_{t\in[0,T]}X(t)>u}=\Ha T \frac{1}{\sqrt{2\pi}}u^{\frac{2}{\alpha}-1}\exp\LT(-\frac{u^2}{2}\RT)\ooo \ \ \ \text{as}\ \ u \rw \IF,
\EQN
where $ \mathcal{H}_{\alpha}$ is the {\it Pickands constant} defined by
\BQNY\label{pick}
\Ha=\lim_{S\rightarrow\infty}\frac{1}{S}\E{ \exp\biggl(\sup_{t\in[0,S]}\Bigl(\sqrt{2}B_\alpha(t)-t^{\alpha}\Bigr)\biggr)} \in (0,\IF), \EQNY
with $\{B_\alpha(t),t\ge0\}$ a standard fractional Brownian motion (fBm) with \rL{Hurst index} $\alpha/2\in (0,1]$. J. Pickands III proved
\re{(1.1)} using the {\it double sum method} and  
the following asymptotics (set $S\in (0,\IF)$)
\BQN\label{eq:pick}
\pk{\sup_{t\in[0,u^{-2/\alpha}S]}X(t)>u}=\Ha[0, S] \frac{1}{\sqrt{2\pi}u}\exp\LT(-\frac{u^2}{2}\RT)\ooo \ \ \ \text{as}\ \ u \rw \IF,
\EQN
where
\BQNY\label{pick}
\Ha[0, S]=\E{ \exp\biggl(\sup_{t\in[0,S]}\Bigl(\sqrt{2}B_\alpha(t)-t^{\alpha}\Bigr)\biggr)}\rL{\in(0,\IF)}.
\EQNY
Piterbarg (1972) obtained a similar result for non-stationary Gaussian processes, namely 
\BQN\label{eq:piter1}
\pk{\sup_{t\in[0,u^{-2/\alpha}S]}\frac{X(t)}{1+dt^\alpha}>u}=\P_{\alpha,\alpha}^{d}[0,S] \frac{1}{\sqrt{2\pi}u}\exp\LT(-\frac{u^2}{2}\RT)\ooo \ \ \ \text{as}\ \ u \rw \IF,
\EQN
 where \rL{$X$ is the centered stationary Gaussian process as above, $d> 0$, and}
\BQNY\label{pick}
\P_{\alpha,\beta}^d[0, S]=\E{ \exp\biggl(\sup_{t\in[0,S]}\Bigl(\sqrt{2}B_\alpha(t)-\abs{t}^{\alpha}-d\abs{t}^{\beta}\Bigr)\biggr)}\in (0,\IF), \ \ \rL{d,\beta,S \in(0,\IF)}.
\EQNY
For a centered non-stationary Gaussian process $\{X(t),t\ge 0\}$
 with a.s. continuous sample paths  the next two assumptions are crucial:

{\bf Assumption A1}.  The standard deviation function $\sigma_X(\cdot)$ of $X$ attains
its maximum \rL{(assumed to be 1)} 
over $[0,T]$ at the unique point $t=T$.  Further, there exist some positive constants
$\nu\in (0,2],\mu,A,D$ such that
\BQN\label{eq:VarX}
\sigma_X(t)&=& 1 -A(T-t)^{\mu}+o((T-t)^{\mu}), \quad t\rw T,
\EQN
and
\BQN\label{eq:rX}
r_X(s,t)=Corr(X(s), X(t))=1- D|t-s|^{\nu}+o(|t-s|^{\nu}), \quad \min(t,s)\rw T.
\EQN
{\bf Assumption A2}. There exist positive constants $ G$ and $\gamma$ such that
\BQN\label{eq:IncX}
\E{(X(t)-X(s))^{2}} &\leq & G|t-s|^{\gamma}
\EQN
holds for all $s,t\in[0, T]$.

For such a centered non-stationary Gaussian process $\{X(t),t\ge 0\}$ it is known that \rL{(see e.g. D\c{e}bicki and Sikora (2011) or Piterbarg (1996))}
\BQNY
\pk{\sup_{t\in[0,T]}X(t)>u}=\DD_{\nu,\mu}  \frac{1}{\sqrt{2\pi}}u^{\LT(\frac{2}{\nu}-\frac{2}{\mu}\RT)_+ -1}\exp\LT(-\frac{u^2}{2}\RT)\ooo \ \ \ \text{as}\ \ u \rw \IF,
\EQNY
where $(x)_+=\max(0,x)$, and $\DD_{\nu,\mu}$ is a positive constant, which, when $\mu=\nu$ is  commonly referred to as
the {\it Piterbarg  constant} defined by
\BQNY\label{piter}
\P_{\nu,\nu}^{\frac{A}{D}}=\lim_{S\rightarrow\infty}\P_{\nu,\nu}^{\frac{A}{D}}[0,S] \in (0,\IF).
\EQNY
\rL{It is worth pointing out that in Theorem D.3 in Piterbarg (1996) it is assumed that the unique maximum point of  $\sigma_X(\cdot)$  is attained at some inner point of $(0,T)$; \re{in that case} the Piterbarg constant is given by $\tilde{\P}_{\nu,\nu}^{\frac{A}{D}}=\lim_{S\rightarrow\infty}\P_{\nu,\nu}^{\frac{A}{D}}[-S,S]$.}

Let  $\{\chi_n(t),t\ge0\}$ be a
chi-process with $n\inn$ degrees of freedom defined by
\BQNY
\chi_n(t)\ =\ \sqrt{\sum_{i=1}^n X_i^2(t)},\quad {t\ge0,}
\EQNY
where $\{X_i(t), t\ge0\}$, $1\le i\le n$, are independent copies of a centered Gaussian process $\{X(t), t\ge0\}$ with a.s. continuous sample paths.
The investigation of
\BQN\label{AlPi}
\pk{\sup_{t\in[0,T]} \chi_n(t) >u}\ \ \ \text{as}\ u \rw \IF
\EQN
was initiated by an envelope of a Gaussian process over a high level, see e.g., Belyaev and Nosko (1969), Lindgren (1980a,b, 1989).
The tail asymptotic behaviour of chi-processes is crucial for numerous statistical applications, see e.g.,  Aronowich and  Adler (1985), Albin and Jaru$\check{s}$kov\'{a} (2003), 
Jaru\v{s}kov\'{a} (2010), Jaru$\check{s}$kov\'{a} and Piterbarg (2011), and the references therein. We mention in passing that the limit behaviour of maximum of chi-processes is the same as that for Gaussian processes (Kabluchko (2011), Hashorva et al.\ (2012));
in the limit the Brown-Resnick process appears.\\
Albin (1990) studied  the exact asymptotics of \eqref{AlPi} for a centered stationary generalized chi-process using Berman's approach (see Berman (1992)), whereas Piterbarg (1994a) obtained a generalization of Albin's  result by resorting to {\it the double sum method}.
In Piterbarg (1994b), the author investigated the exact asymptotics of \eqref{AlPi} for a centered non-stationary generalized chi-process where the generic Gaussian process is differentiable and with variance attaining its global maximum at only one inner point of the interval $[0,T]$. Throughout the paper, a chi-process generated by centered (non-)stationary Gaussian processes is called a {\it (non-)stationary chi-process}.\\

Let $g(\cdot)$ be a non-negative bounded measurable function satisfying one of the following two conditions:

{\bf Assumption G1}.  $g(\cdot)$ attains its minimum  0 over $[0,T]$ at unique point $0$, and further there exist some positive constants $c, \beta$ such that
$$
g(t)=ct^{\beta}\ooo,\ \ \ t\to0;
$$
{\bf Assumption G2}. There exist some constants $ \widetilde{c}\in \R$ and $ \widetilde{\beta}>0$ such that
$$
g(t)=g(T)- \widetilde{c}(T-t)^{ \widetilde{\beta}}\ooo,\ \ \ t\to T.
$$
In this paper, we derive the exact asymptotics of
\BQN\label{eq:mainchi}
\pk{\sup_{t\in[0,T]} \xih >u}\ \ \ \text{as}\  u \rw \IF
\EQN
for
i) stationary chi-processes with a trend function $g(\cdot)$ satisfying Assumption G1;
ii) non-stationary chi-processes a trend function $g(\cdot)$ satisfying Assumption G2.

\rL{The investigation of the tail asymptotics of the maximum of  chi-processes with trend  is motivated by the problem of \re{the} exit of a vector Gaussian load process in engineering sciences, see, e.g., Lindgren (1980a) and the references therein. More precisely, let $\vk{X}(t)=(X_1(t),\cdots,X_n(t)), t\ge0$ be a vector Gaussian load
process. Of interest is the probability of exit
\BQNY
\pk{\vk{X}(t)\not\in\vk{ S}_u(t),\ \ \text{for}\ \text{some}\ t\in[0,T]},
\EQNY
with a {\it time-dependent safety region}
$$
\vk{S}_u(t)=\Biggl\{(x_1,\cdots,x_n)\in\R^n: \sqrt{\sum_{i=1}^n x_i^2}\le h(t,u)\Biggr\}.
$$
The model where $h(t,u)\equiv u$ was considered extensively in the literature as mentioned above; the model where $h(t,u)=u\times c(t)$, with $c(\cdot)$ a positive measurable function, was mentioned in Kozachenko and Moklyachuk (1999) where the authors mainly focused on the exit problem of a class of  square-Gaussian processes. In this paper we shall consider a tractable case that $h(t,u)=u+g(t)$, with $g(\cdot)$ defined as above. The results obtained might also be useful in reliability theory and mathematical statistics applications.}
The analysis of \eqref{eq:mainchi} is based on a tailored {\it double sum method} for chi-processes.
Surprisingly, a {\it generalized Piterbarg constant} $ \mathcal{P}_{\alpha,\beta}^d$, with $\alpha\in(0,2], \beta=\alpha/2, d>0,$ defined by
\BQNY\label{piter}
\P_{\alpha,\beta}^d=\lim_{S\rightarrow\infty}\P^d_{\alpha,\beta}[0,S]\in(0,\IF) 
\EQNY
appears in the asymptotics of the stationary chi-process with trend (we do not observe a generalized Pickands constant as in D\c{e}bicki (2002)). \\ 

Organization of the paper: The main results for the stationary and non-stationary chi-processes with trend are given in Section 2. The proofs are relegated to Section 3 which is followed
then by an Appendix. 

\section{Main Results}
In order to avoid repetitions we shall consider below a chi-process $\{\chi_n(t), t\ge 0\}$ as defined above by taking independent copies of
a generic centered Gaussian processes $X$ with a.s. continuous sample paths. Our asymptotic results will thus depend on the properties of the
Gaussian process $X$. As expected, the stationary case is completely different compared with the non-stationary one.
Throughout this paper denote 
\def\PTT{ \Upsilon_n (u) }
$$ \PTT:= \frac{2^{(2-n)/2}}{\Gamma(n/2)}u^{n-2}\exp\left(-\frac{u^2}{2}\right)$$
which is the asymptotic expansion of the survival function of $\chi_n(0)$ i.e., 
\BQNY
\pk{\chi_n(0)>u }= \PTT \ooo\ \ \ \text{as}\ u \rw \IF,
\EQNY
provided that $X(0)$ is standard normal (i.e., a $N(0,1)$ random variable).

We first present two preliminary results on the tail asymptotics of the maximum of stationary chi-processes  without trend.
The \re{next result} can be found in Piterbarg (1996).
\BS
\label{thP}
Let $\{X(t),t\ge 0\}$ be a stationary Gaussian process with covariance function \rL{$r(\cdot)$}  satisfying assumptions
${\bf R1}-{\bf R2}$ with $\alpha \in (0, 2]$. 
Then, for any constant $T\in(0,\IF)$
\begin{eqnarray}
 \pk{\sup_{t\in[0,T]}\chi_n(t)>u}\ =\ T \mathcal{H}_{\alpha} u^{\frac{2}{\alpha}}\PTT\ooo\label{eqchiu}
\end{eqnarray}
holds as $u\to \infty$.
\ES
\rL{An implication of the last result is the following proposition which will play an important role in the proof of our main results; it can be easily \re{derived by examining} the arguments in Piterbarg (1996).
\BS
\label{thP2}
Let $f(\cdot)$ be a positive function defined in $[0,\IF)$ such that $\lim_{u\to\IF}f(u)/u=1$ and let $S\in(0,\IF)$ be a constant. Under \re{the} assumptions of \neprop{thP} we have that
\begin{eqnarray}
 \pk{\sup_{t\in[0,u^{-2/\alpha}S]}\chi_n(t)>f(u)}\ = \mathcal{H}_{\alpha}[0,S] \Upsilon_n(f(u))\ooo  \label{eqchiu2}
\end{eqnarray}
holds as $u\to \infty$.
\ES
It is worth mentioning that Propositions \ref{thP} and \ref{thP2} are parallel results of Pickands for chi-processes; see \eqref{eq:PP} and \eqref{eq:pick}. }

Next, we give our first result concerning the exact tail asymptotics of the supremum of stationary chi-processes with  trend.

\begin{theo}
\label{MainThm}
Suppose that the covariance function $r(\cdot)$ of the centered stationary Gaussian process $\{X(t),t\ge 0\}$ satisfies assumptions
${\bf R1}-{\bf R2}$ with $\alpha \in (0, 2]$. 
Assume further that $g(\cdot)$ satisfies Assumption G1. If the positive constants $\alpha,\beta,c $ are such that
\BQN\label{eq:c}
c>\left\{
            \begin{array}{ll}
\frac{1}{\beta}, & \hbox{if } \alpha<2\beta ,\\
\frac{2}{\alpha}&  \hbox{if } \alpha\ge2\beta
              \end{array}
            \right.
\EQN
holds,  then
\BQN\label{eq:chiu2}
\pk{\sup_{t\in[0,T]} \xih >u} 
=  \MM_{\alpha,\beta}^c u^{(\frac{2}{\alpha}-\frac{1}{\beta})_+} \PTT
\ooo
\EQN
as $u\to \infty$, where
$$
\MM_{\alpha,\beta}^c=\left\{
            \begin{array}{ll}
c^{-1/\beta}\Gamma(1/\beta+1)\Ha, & \hbox{if } \alpha<2\beta ,\\
\P_{\alpha,\alpha/2}^{c}, & \hbox{if } \alpha=2\beta,\\
1&  \hbox{if } \alpha>2\beta.
              \end{array}
            \right.
$$

\end{theo}

\begin{remarks}
\COM{a) We have the following relation.
\BQNY
\pk{\sup_{t\in[0,T]} \xih >u} =  \MM_{\alpha,\beta}^c u^{(\frac{2}{\alpha}-\frac{1}{\beta})_+} \pk{\chi_n(0)>u }\ooo\ \ \ \text{as}\  u \rw \IF.
\EQNY
}

a) 
For any $d>0$
\BQN
\P_{2,1}^d=\frac{1}{\sqrt{2\pi}}\int_{-\IF}^{\frac{d}{\sqrt{2}}}e^{-\frac{x^2}{2}}dx+\frac{1}{d\sqrt{\pi}}e^{\frac{d^2}{4}-1}.
\EQN
In general $\P_{\alpha,\alpha/2}^d$ is an unknown positive constant which can be eventually calculated by simulations. We mention in passing the paper of Dieker and Yakir (2013) where a new approach is introduced for estimating the Pickands constants.\\
\rL{
b) We see from the proof of last theorem that the minimum of the trend function $g(\cdot)$ taking on $[0,T]$ plays a crucial role. If we assume that $t_0=argmin_{t\in[0,T]}g(t)\in(0,T)$ which is unique and further there exist some positive constants $c, \beta$ such that
$$
g(t)=g(t_0)+c\abs{t}^{\beta}\ooo,\ \ \ t\to t_0,
$$
then \eqref{eq:chiu2} still
holds with $u$ replaced by $u+g(t_0)$, $\Gamma(\cdot)$ replaced by $2\Gamma(\cdot)$, and $\P_{\alpha,\alpha/2}^{c}$ replaced by \\  $\tilde{\P}_{\alpha,\alpha/2}^{c}:=\lim_{S\to\IF} \P_{\alpha,\alpha/2}^{c}[-S,S].$ 
}

c) In view of our proofs and the key results of Piterbarg (1994a) it is possible to obtain additional results for generalized chi-processes. For instance,  if $\{\chi_n(t),t\ge0\}$ is a generalized stationary chi-process  defined by
\BQNY
\chi_n(t)\ =\ \sqrt{\sum_{i=1}^n b_i^2X_i^2(t)},\quad {t\ge0},
\EQNY
with $1=b_1=\cdots=b_k>b_{k+1}\ge b_{k+2}\ge \cdots\ge b_n$, for some $1\le k< n,$ then under assumptions of \netheo{MainThm}
\BQNY
\pk{\sup_{t\in[0,T]} \xih >u} = \prod_{i=k+1}^n(1-b_i^2)^{-1/2} \MM_{\alpha,\beta}^c u^{(\frac{2}{\alpha}-\frac{1}{\beta})_+} \Upsilon_k (u) \ooo
\EQNY
as $u\to \infty$. In order to keep a suitable length of the paper and to avoid extra notation we do not consider here general chi-processes.\\
\end{remarks}

{\bf Examples of $X$}: Numerous important Gaussian processes  satisfy the assumptions of \netheo{MainThm}. We present \re{next}
two interesting cases:

\underline{Fractional Gaussian noise}: Consider $X$ to be the fractional Gaussian noise, i.e.,
$$X(t)=B_\alpha(t+1)-B_\alpha(t), \quad t\ge 0,$$
with $B_\alpha$ a fBm with \rL{Hurst index} $\alpha/2\in(0,1)$. For $\alpha=1$, $X$ is also known as Slepian process.
Clearly $X$ is stationary for any $\alpha \in (0,2)$ and further the covariance function satisfies 
$$ r(t) = 1 - \abs{t}^{\alpha}\ooo,  \quad t \rw 0;\ \ \ \text{and}\ \ r(t)<1\ \text{for all} \ t>0.$$

\underline{Lamperti  transformation of fBm}: Define the Gaussian process $X$ via Lamperti transform of a fBm, i.e., $X(t)=e^{-\alpha/2t}B_\alpha(e^t)$,  which is again a stationary Gaussian process. For the covariance function we have
$$r(t) = 1 - \frac{1}{2}\abs{t}^{\alpha}\ooo,  \quad t \rw 0;\ \ \ \text{and}\ \ r(t)<1\ \text{for all} \ t>0.$$

\COM{
\underline{Ornstein-Uhlenbeck process}: Consider a centered stationary Gaussian process $X$ with correlation function given by
$r(t)=e^{-t}$. We have thus
  $$r(t) = 1 - t \ooo,  \quad t \rw 0;\ \ \ \text{and}\ \ r(t)<1\ \text{for all} \ t>0.$$}

\COM{
\begin{example}
Consider the fractional Gaussian noise $X(t)=B_H(t+1)-B_H(t)$, with $B_H$ being a fBm with \rL{Hurst index} $H\in(0,1]$.  It follows that,  for any $t<1$
$$r(t) = \frac{1}{2}\LT((1+t)^{2H}+(1-t)^{2H}-2t^{2H}\RT).$$
Thus
$$r(t) = 1 - t^{2H} + o(t)\ \ \text{as}\ \  t \rightarrow 0.$$
\end{example}

\begin{example}
Let $X(t)$ be a Ornstein-Uhlenbeck process with covariance function given by
$$\E{X(t)X(0)}=\frac{a}{b}e^{-bt},\ t\ge0,$$
\end{example}
with some constants $a,b>0.$
\begin{example}
Let $X(t)=\int_{t}^{t+1}Z(s)ds$, with $Z(s)$ being a centered stationary Gaussian process.

\end{example}

\begin{example}
Let $X(t)=e^{-Ht}B_H(e^t)$, with $B_H$ being a fBm with \rL{Hurst index} $H\in(0,1]$. It is known that $X(t)$ is a stationary process, e.g., Lamperti (1962) and Samorodnitsky and Taqqu (1994). We see that
\BQNY
r(t)&=& \frac{1}{2}\LT(e^{Ht}+e^{-Ht}-\LT(e^{\frac{1}{2}t}-e^{-\frac{1}{2}t}\RT)^H\RT)\\
&=& 1 - \frac{1}{2}t^H (1+ o(1))\ \ \text{as}\ \  t \rightarrow 0.
\EQNY
\end{example}
}

\bigskip
Next, we deal with a large class of non-stationary chi-processes presenting first the result for chi-process without trend.
\COM{
Let  $\{X_i(t), t\ge0\}$, $1\le i\le n$, be independent copies of a centered Gaussian process $\{X(t), t\ge0\}$ with a.s. continuous sample paths. 
We shall impose two main common assumptions on the Gaussian process $\{X(t), t\ge0\}$.
 Throughout this paper a process $\xi$ with a bar represents a standardized process i.e., $\overline{\xi}(t):= \xi(t)/\sigma_\xi(t)$, with
$\sigma_\xi(\cdot)$ being standard deviation function of the Gaussian process $\{\xi(t),t\ge0\} $.


{\bf Assumption A1}.  The standard deviation function $\sigma_X(\cdot)$  attains
its maximum, denoted by $\widetilde{\sigma}$, over $[0,T]$ at the unique point $t=T$.  Further, there exist some positive constants
$\nu\in (0,2],\mu,A,D$ such that
\BQN\label{eq:VarX}
\sigma_X(t)&=& \wsig -A(T-t)^{\mu}+o((T-t)^{\mu}), \quad t\rw T,
\EQN
and
\BQN\label{eq:rX}
r_X(s,t)=Cov(\overline{X}(s), \overline{X}(t))=1- D|t-s|^{\nu}+o(|t-s|^{\nu}), \quad \min(t,s)\rw T.
\EQN
{\bf Assumption A2}. There exist positive constants $G, \delta<T$ and $\gamma$ such that
\BQN\label{eq:IncX}
\E{(X(t)-X(s))^{2}} &\leq & G|t-s|^{\gamma}
\EQN
holds for for all $s,t\in[T-\delta , T]$.
}
\begin{theo}
\label{MainThm2}
Assume that the centered Gaussian process $\{X(t), t\ge 0\}$ satisfies assumptions {\bf{A1-A2}} for the constants therein. Then for any $T_1\in[0,T)$ we have
\BQN\label{eq:chi2}
\pk{\sup_{t\in[T_1,T]}\chi_n(t)>u}  
&=&  \MM_{\nu,\mu} u^{\LT(\frac{2}{\nu}-\frac{2}{\mu}\RT)_+} \Upsilon_n(u)
\ooo
\EQN
as $u\to \infty$, where
$$
\MM_{\nu,\mu}=\left\{
            \begin{array}{ll}
D^{1/\nu}\frac{\Gamma(1/\mu+1)}{A^{1/\mu}}\mathcal{H}_{\nu}, & \hbox{if } \nu<\mu ,\\
\P_{\nu,\nu}^{\frac{A}{D}}, & \hbox{if } \nu=\mu,\\
1&  \hbox{if } \nu>\mu.
              \end{array}
            \right.
$$

\end{theo}

We state below \re{an extension of Piterbarg  theorem allowing the non-stationary chi-processes to have a non-zero trend.}
\begin{theo}

\label{MainThm3} Assume that $g(\cdot)$ is a positive bounded measurable  function satisfying Assumption G2.
Under the assumptions of \netheo{MainThm2}, 
 if $\mu\le  \widetilde{\beta},$ then  (set $u_*:= u+ g(T)$)
\BQN\label{eq:chi3}
\pk{\sup_{t\in[0,T]} \xih >u} 
=  \rL{\MM_{\nu,\mu}} u_*^{\LT(\frac{2}{\nu}-\frac{2}{\mu}\RT)_+} \Upsilon_n\LT(u_*\RT)
\ooo
\EQN
as $u\to \infty$.
\COM{, where
$$
\WW_{\nu,\mu}=\left\{
            \begin{array}{ll}
D^{1/\nu}\wsig^{\frac{1}{\mu}}\frac{\Gamma(1/\mu+1)}{A^{1/\mu}}\mathcal{H}_{\nu}, & \hbox{if } \nu<\mu ,\\
\DD_{\nu,\mu}, & \hbox{if } \nu\ge \mu.\\
              \end{array}
            \right.
$$
}
\end{theo}

\begin{remarks}
a) As it can be seen from the last two theorems that the only difference between the cases with and without trend is $g(T)$ in $u_*$.
\COM{As it can be seen from our result and its proof, the particular form of the trend function, i.e., $ct^\beta$ is not essential for the exact asymptotics stated in \eqref{eq:chi3}. A more general trend function $g(t)$ can be considered under certain assumptions related to the smoothness of the trend at the point $T$. Our condition $\mu \le 1$ is exactly related to the local properties of $ct^\beta$ at $T$, which in our case do  not depend on $T$.}

b) We conclude from the proof of \netheo{MainThm2} that the Assumption A2 can be relaxed where it can be assumed that there is some $T_0\in(T_1,T)$ such that \eqref{eq:IncX} holds for all $s,t\in[T_0,T]$.
\end{remarks}

{\bf Examples of $X$}: Several important Gaussian processes  satisfy the assumptions of Theorems \ref{MainThm2} and \ref{MainThm3}. We present below
three interesting Gaussian processes (discussed  in Houdr\'{e} and Villa (2003), Bojdecki et al. (2004) and D\c{e}bicki and Tabi\'{s} (2011), respectively).

\underline{Bi-fractional Brownian motion}: Consider $B_{K,H}$ with $K, H\in(0,1)$ to be a bi-fBm, i.e., a self-similar Gaussian process with covariance \re{function} given by
$$
Cov(B_{K,H}(t),B_{K,H}(s))=\frac{1}{2^K}\LT((t^{2H}+s^{2H})^K-\abs{t-s}^{2KH}\RT),\ \ t,s\ge0.
$$
It follows that the standard deviation $\sigma$ of $B_{K,H}$ attaints its maximum over $[0,T]$ at unique point $T$ and
$$
\sigma(t)=T^{KH}-KHT^{KH-1}(T-t)\ooo,   \quad t \rw T.
$$
Further
$$
1-Corr(B_{K,H}(t),B_{K,H}(s))=\frac{1}{2^KT^{2KH}}\abs{t-s}^{2KH}\ooo,    \quad t,s \rw T
$$
and for all $s,t\in[0,T]$ there exists some constant $G>0$ such that
\BQNY
\E{(B_{K,H}(t)-B_{K,H}(s))^{2}} &\leq & G|t-s|^{2KH}.
\EQNY

\underline{Sub-fractional Brownian motion}: The sub-fBm $S_{H}$ with $H\in(0,1)$ is a self-similar Gaussian process with covariance given by
$$
Cov(S_H(t),S_H(s))=t^{2H}+s^{2H}-\frac{1}{2}\LT((s+t)^{2H}-\abs{t-s}^{2H}\RT),\ \ t,s\ge0.
$$
The standard deviation $\sigma$ of $S_H$ attaints its maximum over $[0,T]$ at unique point $T$ and
$$
\sigma(t)=\sqrt{2-2^{2H-1}}T^{H}-\sqrt{2-2^{2H-1}}HT^{H-1}(T-t)\ooo,   \quad t \rw T.
$$
Moreover
$$
1-Corr(S_H(t),S_H(s))=\frac{1}{2(2-2^{2H-1})T^{2H}}\abs{t-s}^{2H}\ooo,    \quad t,s \rw T
$$
and, for all $s,t\in[0,T]$, there exists some constant $G>0$, such that
\BQNY
\E{S_H(t)-S_H(s))^{2}} &\leq & G|t-s|^{H/2}.
\EQNY
\underline{Mean integrated fBm}: Consider a Gaussian process $X_H$ given by
$$
  X_H(t)=\left\{
 \begin{array}{cc}
  {\sqrt{2H+2}\frac{1}{t}\int_{0}^{t}B_{2H}(s)ds},    & t>0,\\
  {0},     & t=0,
 \end{array}
  \right.
$$
with $H\in(0,1)$.
In view of D\c{e}bicki and Tabi\'{s} (2011), we conclude that the standard deviation $\sigma$ of $X_H$ attaints its maximum over $[0,T]$ at unique point $T$ and
$$
\sigma(t)=T^{H}-HT^{H-1}(T-t)\ooo,   \quad t \rw T.
$$
Further
$$
1-Corr(X_H(t),X_H(s))=\frac{1}{2T^{2}}(1-H^2)\abs{t-s}^{2}\ooo,    \quad t,s \rw T
$$
and, for all $s,t\in[\delta,T]$ with some $\delta\in(0,T),$ there exists some constant $G>0$, such that
\BQNY
\E{X_H(t)-X_H(s))^{2}} &\leq & G\delta^{-2}|t-s|.
\EQNY

\COM{
{\bf Example}:
\begin{example}
Let $\{X(t),t\ge 0\}$ be a bi-fractional Brownian motion (bi-fBm) $\{B_{K,H}(t),t\ge 0\}$ with covariance function given by
\BQNY
Cov(B_{K,H}(t),B_{K,H}(s))=\frac{1}{2^{K}}\LT((t^{2H}+s^{2H})^{K}-|s-t|^{2KH}\RT),\quad K\in(0,1], \ H\in(0,1).
\EQNY

\end{example}
}

\def\mbbC{\mathbb{Q}}
\section{Further Results and Proofs}
 In what follows, we give proofs of all the theorems in this paper. Hereafter the positive constant
$\mbbC$ may be different from line to line.

\COM{We first present some preliminary results. The following theorem gives the asymptotics for the survival probability of the supremum of a chi-process without trend, the proof of it can be found in Piterbarg (1996).

\begin{theo}
\label{thP}
Suppose that the covariance function $r(t)$ of the centered stationary Gaussian process $\{X(t),t\ge 0\}$ satisfies assumptions
${\bf R1}-{\bf R2}$ with $\alpha \in (0, 2]$ and some positive $a$. Let $\chi_n(t) =\LT(\sum_{i=1}^n X_i^2(t)\RT)^{1/2}, t\ge0$, where $\{X_i(t),t\ge 0\}, 1\le i\le n$, are independent copies of $\{X(t),t\ge 0\}$. Then
\begin{eqnarray}
\pk{\sup_{t\in[0,T]}\chi_n(t)>u}\ \sim \ Ta^{1/\alpha}\frac{2^{(2-n)/2}\mathcal{H}_{\alpha}}{\Gamma(n/2)}u^{2/\alpha+n-2}\exp\left(-\frac{u^2}{2}\right)\label{eqchiu},\\
\pk{\sup_{t\in[0,u^{-2/\alpha}T]}\chi_n(t)>u}\ \sim \ a^{1/\alpha}\frac{2^{(2-n)/2}}{\Gamma(n/2)}\mathcal{H}_{\alpha}[0,T]u^{n-2}\exp\left(-\frac{u^2}{2}\right)\label{eqchiu2}
\end{eqnarray}
hold as $u\to \infty$.
\end{theo}
Further,
}
Let $\{\xi_u(t,\vkt{v}), t\ge0, \vkt{v}\in\R^{n-1}\}$, $u\ge0$ be a family  of  centered stationary Gaussian \re{random} fields with a.s. continuous sample paths, and covariance function $r_{\xi_u}(t,\vkt{v})$  given by
\BQNY
r_{\xi_u}(t,\vkt{v})=\exp\LT(-u^{-2}D_0t^{\alpha_0}-\sum_{i=1}^{n-1}D_i|v_i|^{\alpha_i}\RT),\ \  t\ge0, \vkt{v}\in\R^{n-1}
\EQNY
for some positive constants $D_i, 0\le i\le n-1,$ and $\alpha_i\in(0,2], 0\le i\le n-1$.


\BT \label{ThmlamS} Let $f(\cdot)$ be a positive function defined in $[0,\IF)$ such that $\lim_{u\to\IF}f(u)/u=1$.
For any $ c, \beta,S_1,S_2>0$ we have
\BQN\label{eq:lamS}
\pk{\underset{\vkt{v}\in \prod_{i=1}^{n-1}\LT[0,u^{-2/\alpha_i}S_2\RT]}{\underset{t\in[0,S_1]}\sup}\frac{\xi_u(t,\vk{v})}{1+ct^\beta u^{-2}}>f(u) }&=& \P_{\alpha_0,\beta}^{cD_0^{-\frac{\beta}{\alpha_0}}}\LT[0,D_0^{\frac{1}{\alpha_0}}S_1\RT]\prod_{i=1}^{n-1}\mathcal{H}_{\alpha_i}\LT[0,D_i^{\frac{1}{\alpha_i}}{S_2}
\RT]\nonumber\\
&&\times\frac{1}{\sqrt{2\pi}f(u)}
\exp\LT(-\frac{(f(u))^2}{2}\RT)(1+\oo(1))
\EQN
 as $u\rw\IF$.
\ET
\prooftheo{ThmlamS} Set $\zeta_u(t,\vkt{v})=\xi_u(t,u^{-2/\alpha_1}v_1,\cdots,u^{-2/\alpha_{n-1}}v_{n-1}), t\ge0, \vkt{v}\in\R^{n-1}, u>0$ with
covariance function 
\BQNY
r_{\zeta_u}(t,\vkt{v})=\exp\LT(-u^{-2}D_0t^{\alpha_0}-u^{-2}\sum_{i=1}^{n-1}D_i|v_i|^{\alpha_i}\RT),\ \  t\ge0, \vkt{v}\in\R^{n-1}, u>0.
\EQNY
Denote further
\BQNY
R_{\zeta_u}(t,\vkt{v},t',\vkt{v}'):=Cov\LT(\frac{\zeta_u(t,\vkt{v})}{1+ct^\beta u^{-2}},\frac{\zeta_u(t',\vkt{v}')}{1+ct'^\beta u^{-2}}\RT) =\frac{r_{\zeta_u}(\abs{t- t'},\vkt{v}-\vkt{v}')}{(1+ct^{\beta}u^{-2})(1+ct'^{\beta}u^{-2})},\ \   t,t'\ge0, \vkt{v},\vkt{v}'\in\R^{n-1}.
\EQNY

Using the classical approach (see e.g., Piterbarg (1996)) we have (set {$\vk{0}=(0 \ldot 0)\inr^{n-1}$})
\BQN\label{eq:lamS1}
&&\pk{\underset{\vkt{v}\in \prod_{i=1}^{n-1}\LT[0,u^{-2/\alpha_i}S_2\RT]}{\underset{t\in[0,S_1]}\sup}\frac{\xi_u(t,\vkt{v})}{1+ct^\beta u^{-2}}>f(u) }
=\frac{1}{\sqrt{2\pi}f(u)}
\exp\LT(-\frac{(f(u))^2}{2}\RT)\nonumber\\
&&\ \ \ \times \int_{-\IF}^\IF e^{w-\frac{w^2}{2(f(u))^2}}\pk{\underset{t\in[0,S_1],\vkt{v}\in \LT[0,{S_2}\RT]^{n-1}}\sup\frac{\zeta_u(t,\vkt{v})}{1+ct^\beta u^{-2}} >f(u) \Bigl| \zeta_u(0,\vk{0})=f(u)-\frac{w}{f(u)}}dw.
\EQN
Further, it follows that
\BQNY
\LT\{\frac{\zeta_u(t,\vkt{v})}{1+ct^\beta u^{-2}} \Bigl | \Bigl(\zeta_u(0,\vk{0}) =f(u)-\frac{w}{f(u)}\Bigr), t\in[0,S_1],\vkt{v}\in \LT[0,{S_2}\RT]^{n-1}\RT\}
\EQNY
has the same distribution as
\BQNY
\LT\{\frac{\zeta_u(t,\vkt{v})}{1+ct^\beta u^{-2}}-R_{\zeta_u}(t,\vkt{v},0,\vk{0}) \zeta_u(0,\vk{0})+R_{\zeta_u}(t,\vkt{v},0,\vk{0})\LT(f(u)-\frac{w}{f(u)}\RT), t\in[0,S_1],\vkt{v}\in \LT[0,{S_2}\RT]^{n-1}\RT\}.
\EQNY
Thus, the integrand in \eqref{eq:lamS1} can be rewritten as
\BQNY
&&\pk{\underset{t\in[0,S_1],\vkt{v}\in \LT[0,{S_2}\RT]^{n-1}}\sup \frac{\zeta_u(t,\vkt{v})}{1+ct^\beta u^{-2}}-R_{\zeta_u}(t,\vkt{v},0,\vk{0}) \zeta_u(0,\vk{0})+R_{\zeta_u}(t,\vkt{v},0,\vk{0})\LT(f(u)-\frac{w}{f(u)}\RT) >f(u)}\\
&&=\pk{\underset{t\in[0,S_1],\vkt{v}\in \LT[0,{S_2}\RT]^{n-1}}\sup \varsigma_u(t,\vkt{v})-(f(u))^2(1-R_{\zeta_u}(t,\vkt{v},0,\vk{0}))+w(1-R_{\zeta_u}(t,\vkt{v},0,\vk{0})) >w},
\EQNY
where
\BQNY
\varsigma_u(t,\vkt{v})=f(u)\LT(\frac{\zeta_u(t,\vkt{v})}{1+ct^\beta u^{-2}}-R_{\zeta_u}(t,\vkt{v},0,\vk{0}) \zeta_u(0,\vk{0})\RT),\quad   t\ge0, \vkt{v}\in\R^{n-1}, u>0.
\EQNY
Next, the following convergence
\BQNY
(f(u))^2(1-R_{\zeta_u}(t,\vkt{v},0,\vk{0}))-w(1-R_{\zeta_u}(t,\vkt{v},0,\vk{0})) \to ct^\beta+D_0t^{\alpha_0}+\sum_{i=1}^{n-1}D_iv_i^{\alpha_i}, \quad u\rw\IF
\EQNY
holds for any $w\in\R$ uniformly with respect to $t\in[0,S_1],\vkt{v}\in\LT[0,{S_2}\RT]^{n-1}$. Furthermore,
\BQNY
\E{\Bigl( \vara_u(t,\vkt{v})-\vara_u(t',\vkt{v}')\Bigr)^2} \to  2D_0|t-t'|^{\alpha_0}+2\sum_{i=1}^{n-1}D_i|v_i-v'_i|^{\alpha_i}, \quad  u\rw\IF
\EQNY
holds uniformly with respect to $t,t'\in[0,S_1],\vkt{v},\vkt{v}'\in\LT[0,{S_2}\RT]^{n-1}$. \rL{ It follows thus that
\BQNY
\E{\Bigl( \vara_u(t,\vkt{v})-\vara_u(t',\vkt{v}')\Bigr)^2} \le  \mathbb{Q}\LT(|t-t'|^{\alpha_0}+ \sum_{i=1}^{n-1} |v_i-v'_i|^{\alpha_i}\RT)
\EQNY
holds for all $u$ sufficiently large and $(t, \vkt{v}),(t',\vkt{v}')$ in any bounded subset of $[0,\IF)\times\R^{n-1}$. Therefore,
the family of random fields $\{\vara_u(t,\vkt{v}),t\in[0,S_1],\vkt{v}\in \LT[0,{S_2}\RT]^{n-1}\}, u>0$ is tight, and thus it  converges weakly to
 $\{\sqrt{2}B_{\alpha_0}(D_0^{1/\alpha_0} t)+\sqrt{2}\sum_{i=1}^{n-1}B_{\alpha_i}(D_i^{1/\alpha_i}v_i),t\in[0,S_1],\vkt{v}\in \LT[0,{S_2}\RT]^{n-1}\}$
 as $u\to\IF$, where $B_{\alpha_i}, i=0,\cdots, n-1$ are independent fBm's with Hurst indexes $\alpha_i/2$, respectively.
 Further using similar arguments as in Lemma 6.1 of Piterbarg (1996) (see also Michna (2009)) we can show that the limit (letting  $u\to\IF$) can be passed under the integral sign in \eqref{eq:lamS1}, and thus the proof is complete.} \QED

Hereafter the diameter of a set $\vk{A}\subset \R^{n}, n\inn$ is defined by
$$
\text{diam}(\vk{A})=\sup_{\vk{t},\vk{s}\in \vk{A}}||\vk{t}-\vk{s}||,
$$
where $||\cdot||$ is the Euclidean norm in $\R^n.$  We write $V_{n}(\vk{A})$ for the $n$-dimensional volume of $\vk{A}$.
\BT \label{ThmSparO} Let $\delta_0$ be a positive constant.
Under the conditions of \netheo{ThmlamS}, for any  $\vk{A}\subset \R^{n-1},n\ge 2,$
with positive volume $V_{n-1}(\vk{A})$
\BQNY
\pk{\underset{t\in[0,S_1],\vkt{v}\in \vk{A}}\sup\frac{\xi_u(t,\vkt{v})}{1+ct^\beta u^{-2}}>u }&=&V_{n-1}(\vk{A})\P_{\alpha_0,\beta}^{cD_0^{-\frac{\beta}{\alpha_0}}}\LT[0,D_0^{\frac{1}{\alpha_0}}S_1\RT]
\prod_{i=1}^{n-1}\mathcal{H}_{\alpha_i}D_i^{\frac{1}{\alpha_i}}\nonumber\\
&&\times\frac{1}{\sqrt{2\pi}} u^{\sum_{i=1}^{n-1}\frac{2}{\alpha_i}-1}\exp\LT(-\frac{u^2}{2}\RT)(1+\oo(1))
\EQNY
holds   as $u\rw\IF$, provided that diam($\vk{A})<\delta_0$ with $\delta_0$ being sufficiently small.
\ET
\rL{\prooftheo{ThmSparO} \re{The proof} follows by similar arguments as in the proof of Lemma 7.1 in Piterbarg (1996) or Lemma 6 in Piterbarg (1994b). \re{It is} mainly based on the double sum method by splitting the set $\vk{A}$ into rectangles and then using Bonferroni's inequality with the aid of \netheo{ThmlamS}. Since it is lengthy and somehow classical, we shall omit the details.\QED}


\def\dU{\delta(u)}
\subsection{Proof of \netheo{MainThm}}
Set next $\dU =\left( \frac{\ln u}{u}\right)^{{1}/{\beta}}, u>0$.  First note that, for any sufficiently small $\vn>0$
\BQNY
\pi_0(u)&:=& \pk{\sup_{t\in\left[\dU,T\right]} \xih >u}\\
&\le& \pk{\sup_{t\in\left[0,T\right]}\chi_n(t)>u+(c-\vn)\frac{\ln u}{u} }
\\
&=&o\LT(u^{n-2+(2/\alpha-1/\beta)_+}\exp\left(-\frac{u^2}{2}\right)\RT)
\EQNY
as $u\rw\IF$, where the last equality follows from \eqref{eqchiu} and the condition \eqref{eq:c}.
Next, we analyze
\BQN\label{eq:lnuu1}
\pk{\sup_{t\in\left[0,\dU\right]} \xih >u},\ \ u\rw\IF,
\EQN
\rL{which, by Assumption G1, is asymptotically equivalent with
\BQN\label{eq:lnuu}
\pi_1(u):=\pk{\sup_{t\in\left[0,\dU\right]} \chicb >u},\ \ u\rw\IF.
\EQN}
It follows from our results below that $\pi_0(u)=o(\pi_1(u))$ as $u\to \IF$. The proof is then established by showing further that $\pi_1(u)$ is asymptotically the same as the right-hand side of \eqref{eq:chiu2}.
\COM{
\BQNY
\pk{\sup_{t\in\left[\dU,T\right]} \xih >u},\ \ u\rw\IF,
\EQNY
is negligible compared with the right-hand side of \eqref{eq:chiu2}.
}
To this end, we need to analyze three cases, namely

i) $\alpha< 2\beta$,\ \ \ ii) $\alpha= 2\beta$,\ \ \ iii) $\alpha> 2\beta$.

\underline{Case  i) $\alpha< 2\beta$}:
Since $\alpha< 2\beta$, for any positive constant $S_1$, we can divide the interval $\LT[0,\dU\RT]$ into several sub-intervals of length $S_1u^{-2/\alpha}$.  Specifically, let for fixed $u>0$
\BQNY
\det_0=u^{-2/\alpha}[0,S_1], \ \ \ \det_k=u^{-2/\alpha}[k S_1,(k +1)S_1], \ \ k \inn.
\EQNY
It follows from Bonferroni's inequality that (set $h(u)=\left\lfloor\frac{(\ln u)^{1/\beta} u^{2/\alpha}}{S_1 u^{1/\beta}} \right\rfloor+1$)
\BQNY
\pi_1(u)&\le& \sum_{k=0}^{h(u)} \pk{\sup_{t\in\det_k} \chicb >u }\\
&\le& \sum_{k=0}^{h(u)} \pk{\sup_{t\in\det_k}\chi_n(t)>u+c(kS_1u^{-2/\alpha})^\beta }\\
&=& \sum_{k=0}^{h(u)} \pk{\sup_{t\in\det_0}\chi_n(t)>u+c(kS_1u^{-2/\alpha})^\beta }=: \pi_2(u).
\EQNY
In view of \eqref{eqchiu2}
\BQN\label{eq:upper1}
\pi_2(u)&=&\frac{2^{(2-n)/2}}{\Gamma(n/2)}\mathcal{H}_{\alpha}[0,S_1]\sum_{k=0}^{h(u)}(u+c(kS_1u^{-2/\alpha})^\beta )^{n-2}\exp\left(-\frac{(u+c(kS_1u^{-2/\alpha})^\beta )^2}{2}\right)(1+ \oo(1))\nonumber\\
&=&\frac{2^{(2-n)/2}}{\Gamma(n/2)}\frac{\mathcal{H}_{\alpha}[0,S_1]}{S_1}u^{2/\alpha-1/\beta+n-2}
\exp\left(-\frac{u^2}{2}\right)\int_0^\IF\exp(-cx^\beta)\, dx (1+ \oo(1))\nonumber\\
&=& \frac{\Gamma(1/\beta+1)}{c^{1/\beta}}\frac{\mathcal{H}_{\alpha}[0,S_1]}{S_1}u^{2/\alpha-1/\beta}
\PTT\ooo
\EQN
as $u\rw\IF$.
Similarly, using  Bonferroni's inequality we obtain
\BQNY
\pi_1(u)&\ge& \sum_{k=0}^{h(u)-1} \pk{\sup_{t\in\det_k} \chicb >u }-\Sigma_{\chi}(u),
\EQNY
where
\BQNY
\Sigma_{\chi}(u):=\sum_{0\le k<j\le h(u)-1} \pk{\sup_{t\in\det_k} \chicb >u, \sup_{t\in\det_j} \chicb >u }.
\EQNY
Along the lines of the proof of \eqref{eq:upper1}, we obtain
\BQN\label{eq:lower1}
\sum_{k=0}^{h(u)-1} \pk{\sup_{t\in\det_k} \chicb >u }
&\ge& \frac{\Gamma(1/\beta+1)}{c^{1/\beta}}\frac{\mathcal{H}_{\alpha}[0,S_1]}{S_1}u^{2/\alpha-1/\beta}
\PTT \ooo
\EQN
as $u\rw\IF$. Furthermore, we have
\BQN\label{eq:doubleS1}
\limsup_{S_1\rw\IF}\limsup_{u\rw\IF}\frac{\Sigma_{\chi}(u)}{A_2(u)}=0, \ \ \ \text{with}\ \  A_2(u):=u^{2/\alpha-1/\beta+n-2}
\exp\left(-\frac{u^2}{2}\right), \ \  u>0.
\EQN
Consequently, the claim for the case  $\alpha< 2\beta$ follows from \eqref{eq:upper1}-\eqref{eq:doubleS1}.
Since the rigorous proof of \eqref{eq:doubleS1} is lengthy, we display it in Appendix.

\underline{Case ii) $\alpha= 2\beta$}:
Clearly,  $S_iu^{-2/\alpha}<\dU$ for $S_i>0, i=1,2,$ when $u$ is sufficiently large. Hence, we have that
\BQNY
\pk{\sup_{t\in[0,S_2u^{-2/\alpha}]} \chicb >u }\le \pi_1(u)
\le \pk{\sup_{t\in\det_0} \chicb >u }+\sum_{k=1}^{h(u)} \pk{\sup_{t\in\det_k} \chicb >u }.
\EQNY
\COM{The following claim that
\BQNY
\pi_1(u) 
&=& \pk{\sup_{t\in\det_0} \chicb >u }\ooo
\EQNY
holds as $u\rw\IF$ and $S_1\rw\IF$ can be established
once we have an expression for
\BQN\label{eq:det0}
\pk{\sup_{t\in\det_0} \chicb >u }, \ \ u\rw\IF
\EQN
and show that
\BQN\label{eq:kIF}
\sum_{k=1}^{h(u)} \pk{\sup_{t\in\det_k} \chicb >u }
=o\LT(\pk{\sup_{t\in\det_0} \chicb >u }\RT)
\EQN
as $u\rw\IF$ and $S_1\rw\IF$.
In order to prove  \eqref{eq:det0},
}
We give  next a technical lemma.
\BEL\label{lemdet} Assume that $\alpha=2\beta$.
 We have
\BQNY
\pk{\sup_{t\in\det_0} \chicb >u }=\P_{\alpha,\beta}^{c}\LT[0,S_1\RT]
\PTT (1+\oo(1))
\EQNY
as $u\rw\IF$.

\EEL
\prooflem{lemdet}
First we see that for $u>0$
\BQN\label{eq:det2}
\pk{\sup_{t\in\det_0} \chicb >u }&=&\pk{\sup_{t\in[0,S_1]} \left(\chi_n(tu^{-2/\alpha})-ct^\beta u^{-1}>u \right)}\nonumber\\
&=&\pk{\sup_{t\in[0,S_1]}\frac{\chi_n(tu^{-2/\alpha})}{1+ct^\beta u^{-2}}>u }.
\EQN
Introduce the Gaussian \re{random} field
$$
Y(t,\vk{s})=\sum_{i=1}^n s_i X_i(t), \quad t\ge0, \vk{s}=(s_1,\cdots,s_n)\in\R^n.
$$
In the light of Piterbarg (1996)
$$
\underset{t\in[0,S_1]}\sup \chi_n(t)=\underset{(t,\vk{s})\in\mathcal{G}_{S_1}}\sup Y(t,\vk{s}),
$$
where $\mathcal{G}_{S_1}=[0,S_1]\times\On$, with $\On$ being the unit sphere (with respect to $L_2$-norm) in $\R^n$.
Therefore, continuing \eqref{eq:det2} for $u>0$  we have
\BQN\label{eq:det3}
\pk{\sup_{t\in\det_0} \chicb >u }
&=&\pk{\underset{(t,\vk{s})\in\mathcal{G}_{S_1}}\sup\eta_u(t,\vk{s})>u },
\EQN
where
$$
\eta_u(t,\vk{s}):= \frac{Y(t u^{-2/\alpha},\vk{s})}{1+ct^\beta u^{-2}}, \ \ t\ge0,\ \vk{s}\in \R^n.
$$
Further, it follows that
\BQNY
Var(\eta_u(t,\vk{s}))=\LT(\frac{1}{1+ct^\beta u^{-2}}\RT)^2, \ \ t\ge0,\ \vk{s}\in \On, u>0
\EQNY
and, for $ t, t'\ge0, \vk{s}, \vk{s'}\in \On$
\BQNY
Corr(\eta_u(t,\vk{s}),\eta_u(t',\vk{s'}))=1-(1-r(u^{-2/\alpha}|t-t'|))-\frac{1}{2}r(u^{-2/\alpha}|t-t'|)||\vk{s}-\vk{s'}||^2.
\EQNY
We split the sphere $\On$ into sets of small diameters $\{\pO_i, 0\le i\le \QQ\}$,  
where
$$
\QQ=\sharp\{\pO_i\}<\IF.
$$
Note that when $n=1$ the sphere $\mathcal{S}_0$ consists of two points $\{1,-1\}$, and thus in this case the partition $\{\pO_i, 0\le i\le 1\}$ consists of two single points. The assertions below is valid for this case as well.
We  have by  Bonferroni's inequality
\BQNY\label{eq:Bon2}
&&\sum_{0\le i\le \QQ}\pk{\underset{t\in[0,S_1],\vk{s}\in\pO_i}\sup\eta_u(t,\vk{s})>u }
\ge\pk{\underset{(t,\vk{s})\in\mathcal{G}_{S_1}}\sup\eta_u(t,\vk{s})>u }\\
&&\quad\ge\sum_{0\le i\le \QQ}\pk{\underset{t\in[0,S_1],\vk{s}\in\pO_i}\sup\eta_u(t,\vk{s})>u }-\sum_{0\le i<l\le \QQ}\pk{\underset{t\in[0,S_1],\vk{s}\in\pO_i}\sup\eta_u(t,\vk{s})>u,
\underset{t\in[0,S_1],\vk{s}\in\pO_l}\sup\eta_u(t,\vk{s})>u }.
\EQNY
We focus next on  $\pO_0$ which includes $(1,0,\cdots,0)$. When diam$(\pO_0)$ is small enough, we can find a one-to-one projection $g$ from $\pO_0$ to the corresponding points where the first component is 1, i.e., $g\vk{v}=(1,v_2,\cdots,v_n)$ for all $\vk{v}=(v_1,v_2,\cdots,v_n)\in\partial O_0$. Thus
\BQNY
\pk{\underset{t\in[0,S_1],\vk{s}\in\pO_0}\sup\eta_u(t,\vk{s})>u }=\pk{\underset{t\in[0,S_1],\vkt{v}\in g\pO_0}\sup\eta_u(t,\vkt{v})>u }.
\EQNY
  Further, in the light of Lemma 10 of Piterbarg (1994b) for any $\vn>0$ there exist positive constants $\delta, u_0$ such that, for diam$(\pO_0)<\delta$, and $u>u_0$
\BQNY
1-\LT(1-\frac{\vn}{2}\RT)  u^{-2}|t-t'|^\alpha-\LT(\frac{1}{2}-\frac{\vn}{4}\RT)\sum_{i=2}^n|s_i-s'_i|^{2}&\ge& Corr(\eta_u(t,\vk{s}),\eta_u(t',\vk{s'}))\\
&\ge&1-\LT(1+\frac{\vn}{2}\RT) u^{-2}|t-t'|^\alpha-\LT(\frac{1}{2}+\frac{\vn}{4}\RT)\sum_{i=2}^n|s_i-s'_i|^{2}
\EQNY
uniformly in $ t, t'\ge0, \vk{s}, \vk{s'}\in \pO_0$.
Define two centered stationary Gaussian processes $\{\xi_u^{\pm}(t,\vkt{v}), t\ge0, \vkt{v}\in \R^{n-1}\}$ with covariance functions given by (set $\vn_{\pm}=1\pm\vn$)
\BQNY
r_{\xi_u^{\pm}}(t,\vkt{v})=\exp\LT(-\vn_\pm  u^{-2}t^\alpha-\frac{\vn_\pm}{2}\sum_{i=1}^{n-1}v_i^{2}\RT),\ \  t\ge0, \vkt{v}\in \R^{n-1},
\EQNY
respectively. In view of Slepian's Lemma (see e.g., Falk et al. (2010)) we have
\BQNY
\pk{\underset{t\in[0,S_1],\vkt{v}\in g\pO_0}\sup\frac{\xi_u^{-}(t,\vkt{v})}{1+ct^\beta u^{-2}}>u }
\le\pk{\underset{t\in[0,S_1],\vkt{v}\in g\pO_0}\sup\eta_u(t,\vkt{v})>u }
\le\pk{\underset{t\in[0,S_1],\vkt{v}\in g\pO_0}\sup\frac{\xi_u^{+}(t,\vkt{v})}{1+ct^\beta u^{-2}}>u }.
\EQNY
Applying \netheo{ThmSparO} to both sides of the last inequality
we conclude that
\BQNY
&&V_{n-1}(g\pO_0)\P_{\alpha,\beta}^{c( \vn_-)^{-\frac{\beta}{\alpha}}}\LT[0,(  \vn_-)^{\frac{1}{\alpha}}S_1\RT]\vn_-^{\frac{n-1}{2}}\frac{1}{(2\pi)^{n/2}}
u^{n-2}\exp\LT(-\frac{u^2}{2}\RT)(1+\oo(1))\le \pk{\underset{t\in[0,S_1],\vkt{v}\in g\pO_0}\sup\eta_u(t,\vkt{v})>u }\\
&&\quad \le V_{n-1}(g\pO_0)\P_{\alpha,\beta}^{c( \vn_+)^{-\frac{\beta}{\alpha}}}\LT[0,(\vn_+)^{\frac{1}{\alpha}}S_1\RT]\vn_+^{\frac{n-1}{2}}\frac{1}{(2\pi)^{n/2}}
u^{n-2}\exp\LT(-\frac{u^2}{2}\RT)(1+\oo(1))
\EQNY
as $u\rw\IF,$ where we used the fact that $\mathcal{H}_2=1/\sqrt{\pi}$. Note that for any sufficiently small positive $\vn_1$, when $\min_{0\le i\le \QQ}\text{diam}(\pO_i)$ is chosen sufficiently small, we have
$$
V_{n-1}(g\pO_i)(1-\vn_1)\le V_{n-1}(\pO_i)\le V_{n-1}(g\pO_i)(1+\vn_1)
$$
for any $0\le i\le \QQ$.
Consequently, by the stationarity of the process $\{\eta_u(t,\vk{s}), (t,\vk{s})\in \mathcal{G}_{S_1}\}$, and then letting $\vn,\vn_1\rw0,$ we conclude that
\BQNY
\sum_{0\le i\le \QQ}\pk{\underset{t\in[0,S_1],\vk{s}\in\pO_i}\sup\eta_u(t,\vk{s})>u }
= V_{n-1}(\On)\P_{\alpha,\beta}^{c}\LT[0,S_1\RT]\frac{1}{(2\pi)^{n/2}}
u^{n-2}\exp\LT(-\frac{u^2}{2}\RT)(1+\oo(1))
\EQNY
as $u\rw\IF$.
Moreover, using similar argumentations as in Appendix we show that
\BQNY
\sum_{0\le i<l\le \QQ}\pk{\underset{t\in[0,S_1],\vk{s}\in\pO_i}\sup\eta_u(t,\vk{s})>u,
\underset{t\in[0,S_1],\vk{s}\in\pO_l}\sup\eta_u(t,\vk{s})>u }=\oo\LT(u^{n-2}\exp\LT(-\frac{u^2}{2}\RT)\RT)
\EQNY
as $u\rw\IF,$ and $S_1\rw\IF$. Since  $ V_{n-1}(\On)=2\pi^{n/2}/\Gamma(n/2)$ the proof is thus complete. \QED

Furthermore, we obtain the following asymptotic upper bound
\BQN\label{eq:kIF1}
\sum_{k=1}^{h(u)} \pk{\sup_{t\in\det_k} \chicb >u }
&\le& \sum_{k=1}^{\IF} \pk{\sup_{t\in\det_k}\chi_n(t)>u+c(kS_1u^{-2/\alpha})^\beta }\nonumber\\
&\overset{\eqref{eqchiu2}}{\le}& \mbbC \ S_1 u^{n-2}\exp\left(-\frac{u^2}{2}\right)\sum_{k=1}^{\IF} e^{-c(kS_1)^{\beta}}\ooo
\EQN
as ${u\rw\IF}$, which together with \nelem{lemdet} yields that, for $S_2>0$
\BQN\label{eq:infsup}
\P_{\alpha,\beta}^{c}\LT[0,S_2\RT]&\le& \liminf_{u\rw\IF}\frac{\pi_1(u)}{\PTT}\nonumber\\
&\le& \limsup_{u\rw\IF}\frac{\pi_1(u)}{\PTT}\le\P_{\alpha,\beta}^{c }\LT[0,S_1\RT]+\mbbC S_1 \sum_{k=1}^{\IF} e^{-c(kS_1)^{\beta}}.
\EQN
Letting $S_2\rw\IF$, we have the finiteness of the {\it generalized Piterbarg constant}, i.e., $\P_{\alpha,\alpha/2}^{c }<\IF$. Similarly, \rL{ letting $S_1\rw\IF$ we obtain $\P_{\alpha,\alpha/2}^{c }>0$.} Consequently,
the claim for the case $\alpha=2\beta$ follows by letting $S_1,S_2\rw\IF$.

\COM{
In order to analyze the two bounds of the last formula, we give next a lemma.
\BEL \label{lemlamS}
For any $S_2, S>0,$ we have, as $u\rw\IF,$
\BQN\label{eq:lamS}
\pk{\underset{t\in[0,S],\vk{s}\in \LT[-\frac{S_2}{u},\frac{S_2}{u}\RT]^{n-1}}\sup\frac{\xi_u^{-}(t,\vk{s})}{1+ct^\beta u^{-2}}>u }= \P_{B_\alpha}^{ct^\beta}[0,S]\LT(\mathcal{H}_2\LT[-\frac{S_2}{\sqrt{2}},\frac{S_2}{\sqrt{2}}\RT]\RT)^{n-1}\frac{1}{\sqrt{2\pi}u}
\exp\LT(-\frac{u^2}{2}\RT)(1+\oo(1)).
\EQN
\EEL
\prooflem{lemlamS} Set $\zeta_u(t,\vk{s})=\xi_u^{-}(t,u^{-1}\vk{s})$. Thus
\BQNY
r_{\zeta_u}(t,\vk{s},t',\vk{s'})=\exp\LT(1-(1\pm\vn)u^{-2}|t-t'|^\alpha-\LT(\frac{1}{2}\pm\frac{\vn}{2}\RT)u^{-2}\sum_{i=2}^n|s_i-s'_i|^{2}\RT).
\EQNY
Denote
\BQNY
R_{\zeta_u}(t,\vk{s},t',\vk{s'}):=Cov\LT(\frac{\zeta_u(t,\vk{s})}{1+ct^\beta u^{-2}},\frac{\zeta_u(t',\vk{s'})}{1+ct'^\beta u^{-2}}\RT) =\frac{r_{\zeta_u}(t,\vk{s},t',\vk{s'})}{(1+ct^{\beta}u^{-2})(1+ct'^{\beta}u^{-2})}.
\EQNY

Utilising the classical approach (see e.g., Piterbarg (1996)) we have for $u\ge0$ 
\BQN\label{eq:lamS1}
&&\pk{\underset{t\in[0,S],\vk{s}\in \LT[-\frac{S_2}{u},\frac{S_2}{u}\RT]^{n-1}}\sup\frac{\xi_u^{-}(t,\vk{s})}{1+ct^\beta u^{-2}}>u }\nonumber\\
&=&\frac{1}{\sqrt{2\pi}u}
\exp\LT(-\frac{u^2}{2}\RT)\int_{-\IF}^\IF e^{w-\frac{w^2}{2u^2}}\pk{\underset{t\in[0,S],\vk{s}\in \LT[-{S_2},{S_2}\RT]^{n-1}}\sup\frac{\zeta_u(t,\vk{s})}{1+ct^\beta u^{-2}} >u| \zeta_u(0,\vk{0})=u-\frac{w}{u}}dw.
\EQN
It can be shown that
\BQNY
\LT\{\frac{\zeta_u(t,\vk{s})}{1+ct^\beta u^{-2}}| \zeta_u(0,\vk{0})=u-\frac{w}{u}; t\in[0,S],\vk{s}\in \LT[-{S_2},{S_2}\RT]^{n-1}\RT\}
\EQNY
has the same distribution as
\BQNY
\LT\{\frac{\zeta_u(t,\vk{s})}{1+ct^\beta u^{-2}}-R_{\zeta_u}(t,\vk{s},0,\vk{0}) \zeta_u(0,\vk{0})+R_{\zeta_u}(t,\vk{s},0,\vk{0})\LT(u-\frac{w}{u}\RT); t\in[0,S],\vk{s}\in \LT[-{S_2},{S_2}\RT]^{n-1}\RT\}.
\EQNY
The integral in \eqref{eq:lamS1} can be rewritten as
\BQNY
&&\pk{\underset{t\in[0,S],\vk{s}\in \LT[-{S_2},{S_2}\RT]^{n-1}}\sup \frac{\zeta_u(t,\vk{s})}{1+ct^\beta u^{-2}}-R_{\zeta_u}(t,\vk{s},0,\vk{0}) \zeta_u(0,\vk{0})+R_{\zeta_u}(t,\vk{s},0,\vk{0})\LT(u-\frac{w}{u}\RT) >u}\\
&&=\pk{\underset{t\in[0,S],\vk{s}\in \LT[-{S_2},{S_2}\RT]^{n-1}}\sup \varsigma_u(t,\vk{s})-u^2(1-R_{\zeta_u}(t,\vk{s},0,\vk{0}))+w(1-R_{\zeta_u}(t,\vk{s},0,\vk{0})) >w},
\EQNY
where
\BQNY
\varsigma_u(t,\vk{s})=u\LT(\frac{\zeta_u(t,\vk{s})}{1+ct^\beta u^{-2}}-R_{\zeta_u}(t,\vk{s},0,\vk{0}) \zeta_u(0,\vk{0})\RT).
\EQNY
It is derived that
\BQNY
\lim_{u\rw\IF}u^2(1-R_{\zeta_u}(t,\vk{s},0,\vk{0}))-w(1-R_{\zeta_u}(t,\vk{s},0,\vk{0}))=ct^\beta+ (1-\vn)t^\alpha-\LT(\frac{1}{2}-\frac{\vn}{2}\RT)\sum_{i=2}^ns_i^{2},
\EQNY
uniformly on $w\in\R$, and
\BQNY
\lim_{u\rw\IF}\E{\vara_u(t,\vk{s})-\vara_u(t',\vk{s'})}^2= 2(1-\vn)|t-t'|^\alpha-\LT(1-{\vn}\RT)\sum_{i=2}^n|s_i-s'_i|^{2},
\EQNY
uniformly on $t,t'\in[0,S],\vk{s},\vk{s'}\in\LT[-{S_2},{S_2}\RT]^{n-1}$.
Therefore, we conclude that the process $\{\vara_u(t,\vk{s}), t\in[0,S],\vk{s}\in \LT[-{S_2},{S_2}\RT]^{n-1}\}$ converge weakly to the process $\{\sqrt{2(1-\vn)}B_\alpha(t)+\sum_{i=2}^n\NN_i s_i\}, t\in[0,S],\vk{s}\in \LT[-{S_2},{S_2}\RT]^{n-1}$, where $B_\alpha(t)$ is the fBm with \rL{Hurst index} $\alpha/2\in(0,1]$, and $\NN_i, i=2,\cdots,n,$ are independent standard normal distributions, which are independent of $B_\alpha$.

Consequently, the claim follows from the above discussions. \QED

\BEL \label{lemSparO}
When diam($g\parO)<\delta$, we have, as $u\rw\IF,$
\BQNY
\pk{\underset{t\in[0,S],\vk{s}\in g\parO}\sup\frac{\xi_u^{-}(t,\vk{s})}{1+ct^\beta u^{-2}}>u }=V_{n-1}(g\parO)\P_{B_\alpha}^{ct^\beta}[0,S]\frac{1}{(2\pi)^{n/2}}u^{n-2}\exp\LT(-\frac{u^2}{2}\RT)(1+\oo(1)),
\EQNY
where $V_{n-1}(g\parO)$ is the $n-1$-dimensional volume of $g\parO$.
\EEL
\prooflem{lemSparO} Let
$$ \wtre_0=\LT[-\frac{S_2}{u},\frac{S_2}{u}\RT]^{n-1},\  \wtre_{\vk{k}}=\prod_{i=1}^{n-1}\LT[k_i\frac{S_2}{u},(k_i+1)\frac{S_2}{u}\RT],\  \vk{k}\in \mathbb{Z}^{n-1}.$$
Let
$$\CK^{-}=\{\vk{k}:\wtre_{\vk{k}}\in g\parO\},\ N_u^{-}=\sharp\{\CK^{-}\}; \ \CK^{+}=\{\vk{k}:\wtre_{\vk{k}}\cap g\parO\neq \emptyset\}, \ N_u^{+}=\sharp\{\CK^{+}\}.$$
It thus follows that
\BQN\label{eq:NV}
N_u^{-}=N_u^{+}(1+\oo(1))=V_{n-1}(g\parO)\LT(\frac{u}{2S_2}\RT)^{n-1}(1+\oo(1)), \ \ u\rw\IF.
\EQN
Furthermore, Bonferroni's inequality gives that
\BQNY
\sum_{\vk{k}\in \CK^{+}}\pk{\underset{t\in[0,S],\vk{s}\in \wtre_{\vk{k}}}\sup\frac{\xi_u^{-}(t,\vk{s})}{1+ct^\beta u^{-2}}>u }&\ge&\pk{\underset{t\in[0,S],\vk{s}\in g\parO}\sup\frac{\xi_u^{-}(t,\vk{s})}{1+ct^\beta u^{-2}}>u }\\
&\ge&\sum_{\vk{k}\in \CK^{-}}\pk{\underset{t\in[0,S],\vk{s}\in \wtre_{\vk{k}}}\sup\frac{\xi_u^{-}(t,\vk{s})}{1+ct^\beta u^{-2}}>u }\\
&&-
\sum_{\vk{k},\vk{j}\in \CK^{-}}\pk{\underset{t\in[0,S],\vk{s}\in \wtre_{\vk{k}}}\sup\frac{\xi_u^{-}(t,\vk{s})}{1+ct^\beta u^{-2}}>u, \underset{t\in[0,S],\vk{s}\in \wtre_{\vk{j}}}\sup\frac{\xi_u^{-}(t,\vk{s})}{1+ct^\beta u^{-2}}>u }
\EQNY
In view of \eqref{eq:lamS} and \eqref{eq:NV}, we conclude that, as $u\rw\IF,$
\BQN\label{eq:upper2}
&&\sum_{\vk{k}\in \CK^{+}}\pk{\underset{t\in[0,S],\vk{s}\in \wtre_{\vk{k}}}\sup\frac{\xi_u^{-}(t,\vk{s})}{1+ct^\beta u^{-2}}>u }=
N_u^+ \pk{\underset{t\in[0,S],\vk{s}\in \wtre_{0}}\sup\frac{\xi_u^{-}(t,\vk{s})}{1+ct^\beta u^{-2}}>u }\\
&=&V_{n-1}(g\parO)\LT(\frac{u}{2S_2}\RT)^{n-1} \P_{B_\alpha}^{ct^\beta}[0,S]\LT(\mathcal{H}_2\LT[-\frac{S_2}{\sqrt{2}},\frac{S_2}{\sqrt{2}}\RT]\RT)^{n-1}\frac{1}{\sqrt{2\pi}u}
\exp\LT(-\frac{u^2}{2}\RT) (1+\oo(1)).
\EQN
In addition, it is known that
$$
\lim_{S_2\rw\IF}\frac{1}{\sqrt{2}S_2}\mathcal{H}_2\LT[-\frac{S_2}{\sqrt{2}},\frac{S_2}{\sqrt{2}}\RT]=\frac{1}{\sqrt{\pi}}.
$$
Therefore, letting $S_2\rw\IF,$ in \eqref{eq:upper2}, we obtain the asymptotic upper bound
\BQNY
\sum_{\vk{k}\in \CK^{+}}\pk{\underset{t\in[0,S],\vk{s}\in \wtre_{\vk{k}}}\sup\frac{\xi_u^{-}(t,\vk{s})}{1+ct^\beta u^{-2}}>u }
=V_{n-1}(g\parO)\P_{B_\alpha}^{ct^\beta}[0,S]\frac{1}{(2\pi)^{n/2}}u^{n-2}\exp\LT(-\frac{u^2}{2}\RT)(1+\oo(1)).
\EQNY

Finally, we need to analyze the double sum part. \QED
}

\underline{Case iii) $\alpha> 2\beta$}:
The lower bound follows immediately since
\BQNY
\pi_1(u)
&\ge& \pk{\chi_n(0)>u }=\PTT.
\EQNY
In view of \nelem{lemdet} we derive an upper bound as follows 
\BQNY
\limsup_{u\rw\IF}\frac{\pk{\sup_{t\in\left[0,\dU \right]} \chicb >u}}{\PTT}
&\le& \limsup_{u\rw\IF}\frac{\pk{\sup_{t\in\det_0} \chicb >u}}{\PTT}= \rL{\P_{\alpha,\frac{\alpha}{2}}^{c}}\LT[0,S_1\RT].
\EQNY
The proof is completed by letting $S_1\rw0.$

\subsection{Proof of \netheo{MainThm2}}
In this subsection, we give the skeleton of the proof of \netheo{MainThm2} which is based on the double sum method.  Again, we introduce a Gaussian random field
$$
Y(t,\vk{s})=\sum_{i=1}^n s_i X_i(t), \quad t\ge0, \vk{s}=(s_1,\cdots,s_n)\in\R^n.
$$
\re{Since}
$$
\underset{t\in[T_1,T]}\sup \chi_n(t)=\underset{(t,\vk{s})\in[T_1,T]\times\On}\sup Y(t,\vk{s})
$$ for any $T_1\in[0,T)$.
For $t,s\ge0, \vk{v},\vk{w}\in \On$
$$
Var(Y(t,\vk{v}))=\sigma_X^2(t),\ \ \text{and}\ Corr(Y(t,\vk{v}),Y(s,\vk{w}))=r_X(s,t)-\frac{1}{2}r_X(s,t)||\vk{v}-\vk{w}||^2.
$$
 Consequently, by \eqref{eq:VarX}--\eqref{eq:IncX} there is some $\delta\in(0,T)$ close to $T$ such that
  \BQNY
Var(Y(t,\vk{v}))\le  1- A q^\mu(u), \quad \text{with}\ \re{q(u) =\LT(\frac{\ln u}{u}\RT)^{2/\mu}} 
\EQNY
holds for all $t\in \LT[\delta,T-  q(u)  \RT]$
and $\vk{v}\in \On$ when $u$ is sufficiently large, and further for $t,s\in [\delta,T]$ and $\vk{v},\vk{w}\in \On$
\BQNY
\E{(Y(t,\vk{v})-Y(s,\vk{w}))^2}\le \mbbC \LT(|t-s|^\gamma+||\vk{v}-\vk{w}||^2\RT).
\EQNY
Therefore, by Piterbarg inequality (cf. Theorem 8.1 of Piterbarg (1996) or Theorem 8.1 in the seminal paper Piterbarg (2001))
\BQN\label{eq:T0}
\Pi_1(u):=\pk{\underset{t\in \LT[\delta,T-  q(u)  \RT]}\sup \chi_n(t)>u}\le \mbbC\ u^{2/\gamma+n-1}\exp\LT(\frac{u^2}{2\LT(1-A q^\mu(u)\RT)}\RT).
\EQN
Furthermore, we have from Borell-TIS inequality (e.g., Adler and Taylor (2007))
\BQN
\Pi_2(u):=\pk{\underset{t\in \LT[T_1,\delta  \RT]}\sup \chi_n(t)>u}&\le& 
\pk{\sup_{(t,\vk{v})\in \mathcal{G}_{\delta}}Y(t,\vk{v})>u}\nonumber\\
&\le&\exp\LT(-\frac{(u-\mathbb{C})^2}{2\sigma_X^2(\delta)}\RT).
\EQN
Next, we focus on the asymptotics of
\def\wchi{\widetilde{\chi}_n}
\BQNY
\Pi_3(u):=\pk{\underset{t\in \LT[T-  q(u)  ,T\RT]}\sup \chi_n(t)>u}=\pk{\underset{t\in \LT[0,  q(u)  \RT]}\sup \wchi(t)>u},\ \ u\rw\IF,
\EQNY
where $\wchi(t)=\chi_n(T-t),$ for $t\in \LT[0,  q(u)  \RT].$ \rL{From the 
results below we conclude that
\BQN\label{eq:pi123}
\Pi_1(u)=o\LT(\Pi_3(u)\RT),\ \ \Pi_2(u)=o\LT(\Pi_3(u)\RT)
\EQN
as $u\to\IF.$ The proof is thus established by showing further that $\Pi_3(u)$ is asymptotically the same as the right-hand side of \eqref{eq:chi2}.}

Similar to the proof of \netheo{MainThm} we need to distinguish between the following three cases:

i) $\nu<\mu$,\ \ \ ii) $\nu=\mu$, \ \ \ iii) $\nu>\mu$.

\def\tUU{\theta(u)}
Let, for $S_1>0$
$$
\det_0=u^{-2/\nu}[0,S_1],\ \ \ \det_k=u^{-2/\nu}[kS_1,(k+1)S_1],\ k\in\N,
$$
and define $\theta(u)=\LT\lfloor\frac{(\ln u)^{2/\mu} u^{2/\nu}}{S_1 u^{ 2/\mu}}\RT\rfloor+1$.\\

\underline{Case i) $\nu<\mu$}:
Since $\nu<\mu$, using Bonferroni's inequality, we have
\BQNY
\sum_{k=0}^{ \tUU } \pk{\sup_{(t,\vk{v})\in\det_k\times\On} Z(t,\vk{v}) >u }&=&\sum_{k=0}^{ \tUU } \pk{\sup_{t\in\det_k}\wchi(t)>u }\\
&\ge& \pk{\underset{t\in \LT[0,  q(u)  \RT]}\sup \wchi(t)>u}\\
&\ge& \sum_{k=0}^{ \tUU -1  } \pk{\sup_{t\in\det_k}\wchi(t)>u }-\Sigma_{\wchi}(u),
\EQNY
where $ Z(t,\vk{v})= Y(T-t,\vk{v}),$ for $ (t,\vk{v})\in \LT[0,  q(u)  \RT]\times \On$, and
\BQNY
\Sigma_{\wchi}(u):=\sum_{0\le k<j\le  \tUU -1  } \pk{\sup_{t\in\det_k}\wchi(t)>u, \sup_{t\in\det_j}\wchi(t)>u }.
\EQNY
\COM{The sum on the left-hand side of the last formula can be rewritten as
\BQNY
\sum_{k=0}^{ \tUU } \pk{\sup_{t\in\det_k}\wchi(t)>u }
&=&\sum_{k=0}^{ \tUU } \pk{\sup_{(t,\vk{v})\in\det_k\times\On} Z(t,\vk{v}) >u },
\EQNY
where}
For any  $\vn \in (0,1)$, when $u$ is sufficiently large, we have
$$
1-A(1-\vn)t^\mu>Var(Z(t,\vk{v}))^{1/2}>1-A(1+\vn)t^\mu
$$
and
\def\sWW{ ||\vk{v}- \vk{w}||^2}
\BQNY
2 D(1-\vn)|t-s|^\nu+ (1-\vn)\sWW &\le&\E{(Z(t,\vk{v})-Z(s,\vk{w}))^2}\\
&\le &2 D(1+\vn)|t-s|^\nu+ (1+\vn)\sWW .
\EQNY
Next we introduce a centered stationary Gaussian process $\{\xi(t),t\ge0\}$ with covariance function
$$r_\xi(t)=\exp\LT(-Dt^\nu\RT),\ \ t\ge0$$
and set
$$Z_2(t,\vk{v})=\sum_{i=1}^n v_i \xi_i(t), \ \ t\ge0,\ \vk{v}\in \R^n,$$
with $\{\xi_i(t),t\ge0\}, 1\le i\le n$, being independent copies of $\{\xi(t),t\ge0\}$. Thus, we have, for $ (t,\vk{v})\in \LT[0,  q(u)  \RT]\times \On$, and $u$  sufficiently large
\BQNY
2D(1-\vn)|t-s|^\nu+(1-\vn)\sWW &\le&\E{(Z_2(t,\vk{v})-Z_2(s,\vk{w}))^2}\\
&\le &2D(1+\vn)|t-s|^\nu+(1+\vn)\sWW .
\EQNY
Therefore, by the arbitrariness of $\vn$ and with the aid of Slepian's Lemma, we conclude that
\BQN\label{eq:ZZ2}
\pk{\sup_{(t,\vk{v})\in\det_k\times\On} Z(t,\vk{v}) >u }=\pk{\sup_{(t,\vk{v})\in\det_k\times\On} Z_2(t,\vk{v})(1-At^\mu) >u }(1+ o(1))
\EQN
as $u\rw\IF$. Consequently, as $u\rw\IF$
\BQNY
\sum_{k=0}^{ \tUU } \pk{\sup_{t\in\det_k}\wchi(t)>u }
&\le& \sum_{k=0}^{ \tUU }\pk{\sup_{(t,\vk{v})\in\det_k\times\On}Z_2(t,\vk{v}) >\frac{u }{1-A(kS_1u^{-2/\nu})^\mu} }\ooo=: \pi_3(u).
\EQNY
Utilising further \eqref{eqchiu2}, we obtain
\BQN\label{eq:upper2}
\pi_3(u)&=&\frac{2^{(2-n)/2}}{\Gamma(n/2)}\mathcal{H}_{\nu}[0,\rL{D^{1/\nu}S_1}]\sum_{k=0}^{ \tUU }\LT(\frac{u}{1-A(kS_1u^{-2/\nu})^\mu} \RT)^{n-2}\exp\left(-\frac{u^2\LT(1+ A (kS_1u^{-2/\nu})^\mu \RT)^2}{2 }\right)(1+ o(1))\nonumber\\
&=&\frac{\mathcal{H}_{\nu}[0,D^{1/\nu} S_1]}{S_1}u^{2/\nu-2/\mu} \Upsilon_n\left( u  \right) \int_0^\IF\exp\LT(- A x^\mu\RT)\, dx (1+o(1))\nonumber\\
&=& D^{1/\nu} \frac{\Gamma(1/\mu+1)}{A^{1/\mu}}\frac{\mathcal{H}_{\nu}[0,D^{1/\nu} S_1]}{D^{1/\nu} S_1} u^{2/\nu-2/\mu}
\Upsilon_n\left( u  \right)(1+o(1))
\EQN
as $u\rw\IF$.
Using the same argumentations as \eqref{eq:upper2} the following asymptotic lower bound
\BQN\label{eq:lower2}
\sum_{k=0}^{ \tUU -1  } \pk{\sup_{t\in\det_k}\wchi(t)>u }
\ge D^{1/\nu} \frac{\Gamma(1/\mu+1)}{A^{1/\mu}}\frac{\mathcal{H}_{\nu}[0,D^{1/\nu} S_1]}{D^{1/\nu} S_1} u^{2/\nu-2/\mu}
\Upsilon_n\left( u  \right)(1+o(1))
\EQN
holds as $u\rw\IF$. 
Furthermore, we have
\BQN\label{eq:doubleS2}
\limsup_{S_1\rw\IF}\limsup_{u\rw\IF}\frac{\Sigma_{\wchi}(u)}{A_3(u)}=0, \ \ \ \text{with}\ \  A_3(u):=u^{2/\nu-2/\mu+n-2}
\exp\left(-\frac{u^2}{2}\right), \ \  u>0
\EQN
the proof of
which is omitted since it is similar to \eqref{eq:doubleS1}. Consequently, the claim for the case $\nu<\mu$ follows from \eqref{eq:upper2}-\eqref{eq:doubleS2}.

\underline{Case ii) $\nu=\mu$}:
Since  $S_iu^{-2/\nu}<  q(u) =\LT(\frac{\ln u}{u}\RT)^{2/\mu} $ for $S_i>0, i=1,2,$ when $u$ is sufficiently large. Hence, we have that
\BQNY
\pk{\sup_{t\in[0,S_2u^{-2/\nu}]} \wchi(t) >u }\le \pi_1(u)
\le \pk{\sup_{t\in\det_0} \wchi(t) >u }+\sum_{k=1}^{\theta(u)} \pk{\sup_{t\in\det_k} \wchi(t) >u }.
\EQNY
\COM{
Since  $S_1u^{-2/\nu}<  q(u)  $ when $u$ is sufficiently large, we have
\BQNY
\pk{\sup_{t\in\left[0,\left(\frac{\ln u}{u}\right)^{{2}/{\mu}}\right]}\wchi(t)>u}&\sim& \pk{\sup_{t\in\det_0}\wchi(t)>u } 
\EQNY
as $u\rw\IF$ and  $S_1\rw\IF$,
once we have the expression of
\BQN\label{eq:det02}
\pk{\sup_{t\in\det_0}\wchi(t)>u }, \ \ u\rw\IF,
\EQN
and show that
\BQN\label{eq:kIF2}
\sum_{k=1}^{ \tUU } \pk{\sup_{t\in\det_k}\wchi(t)>u }
=o\LT(\pk{\sup_{t\in\det_0}\wchi(t)>u }\RT) 
\EQN
as $u\rw\IF$ and  $S_1\rw\IF$.}
From  \eqref{eq:ZZ2}  we obtain further
\BQNY
\pk{\sup_{t\in\det_0}\wchi(t)>u }&=& \pk{\sup_{(t,\vk{v})\in\det_0\times\On} Z_2(t,\vk{v})(1-At^\mu) >u }(1+o(1))\nonumber\\
&=& \pk{\sup_{(t,\vk{v})\in\det_0\times\On} \frac{Z_2(t,\vk{v})}{(1+ A t^\mu)} > u  } (1+o(1))
\EQNY
as $u\rw\IF$. In view of \netheo{ThmlamS} and \netheo{ThmSparO}, and the derivation of the case $\alpha=2\beta$ in the last subsection, we conclude that
\BQN\label{eq:det00}
\pk{\sup_{t\in\det_0}\wchi(t)>u }&=& \P_{\nu,\mu}^{ A D^{-1}}[0,D^{\frac{1}{\nu}}S_1 ]
\Upsilon_n \LT( u \RT)(1+o(1))
\EQN
as $u\rw\IF$. Now, the claim follows using the same argumentation as \eqref{eq:infsup}.

\underline{Case iii) $\nu>\mu$:}
By \eqref{eq:det00} the upper bound is derived as
\BQNY
\limsup_{u\rw\IF}\frac{\pk{\sup_{t\in\left[0,q(u)\right]}\wchi(t)>u}}{\Upsilon_n \LT( u \RT)}
&\le& \limsup_{u\rw\IF}\frac{\pk{\sup_{t\in\det_0}\wchi(t)>u}}{\Upsilon_n \LT( u \RT)}\le \P_{\nu,\rL{\nu}}^{ A D^{-1}}[0,D^{\frac{1}{\nu}}S_1 ].
\EQNY
Since further
\BQNY
\pk{\sup_{t\in\left[0,q(u)\right]}\wchi(t)>u}&\ge& \pk{\wchi(0)>u }=\Upsilon_n \LT( u \RT)
\EQNY
the proof of this case is established by letting $S_1\rw0.$ Consequently, it follows that \eqref{eq:pi123}
\COM{that
\BQNY
\pk{\underset{t\in \LT[\delta,T-  q(u)  \RT]}\sup \chi_n(t)>u}=\oo\LT(u^{n-2+(2/\nu-2/\mu)_+}
\exp\left(-\frac{u^2}{2 }\right)\RT)
\EQNY
as $u\rw\IF$,}
is \re{valid}, and thus the proof is  complete.\QED

\subsection{Proof of \netheo{MainThm3}}
For $\delta\in(0,T)$, set
\BQNY
\Pi(u):=\pk{\sup_{t\in[\delta,T]} \xih >u}.
\EQNY
{Thus for any $u\ge0$
\BQN\label{eq:TT}
\Pi(u) \le\pk{\sup_{t\in[0,T]} \xih >u}
\le\pk{\sup_{t\in[0,\delta]} \xih >u}+\Pi(u).
\EQN
}
It  follows that
\BQNY
\Pi(u)&=&\pk{\sup_{t\in[\delta,T]}\frac{\chi_n(t)}{u+g(t) }>1}\\
&=&\pk{\sup_{t\in[\delta,T]}\frac{\chi_n(t)}{\sigma_X(t)}\frac{m_u(T)}{m_u(t)}>m_u(T)},
\ \ \text{with}\
m_u(t):=\frac{u+g(t)}{\sigma_X(t)}, \  t\ge0.
\EQNY
For any $t \in [0, T]$
$$1-\frac{m_{u}(T)}{m_{u}(t)}=\frac{\sigma_{X}(T)-\sigma_{X}(t)}
{\sigma_{X}(T)}+\frac{\sigma_{X}(t) (g(t)-g(T))}{(u+ g(t) )\sigma_{X}(T)}.$$
Further, in view of \eqref{eq:VarX} and Assumption G2, and noting that $\mu\le \widetilde{\beta}$, $\delta$ can be chosen close enough to $T$ such that
\BQNY
\abs{g(T)-g(t)}\le \mbbC\ (\sigma_{X}(T)-\sigma_{X}(t))
\EQNY
for all $t\in[\delta, T]$. Hence for any  $\vn\in(0,1)$, when $u$ is sufficiently large, we have, uniformly in $[\delta,T]$
\begin{eqnarray}
\label{eq3.1}
1-(1+ \varepsilon)\frac{\sigma_{X}(T)-\sigma_{X}(t)}{\sigma_{X}(T)}
\leq\frac{m_{u}(T)}{m_{u}(t)}
\leq 1-(1-\varepsilon)\frac{\sigma_{X}(T)-\sigma_{X}(t)}{\sigma_{X}(T)}.
\end{eqnarray}
Therefore, for $u$ sufficiently large
$$\pi_{+\vn}(u):=\pk{\sup_{t\in[\delta,T]}Y_{+\vn}(t)>m_{u}(T)}\le \Pi(u)\leq \pi_{-\vn}(u):=\pk{\sup_{t\in[\delta,T]}Y_{-\vn }(t)>m_{u}(T)},$$
where
$$
Y_{\pm \vn }(t)=\sqrt{\sum_{i=1}^n Y_{\pm \vn,i }^2(t)},\ \ t\ge0,
$$
with
$$Y_{\pm \vn, i }(t):=\frac{X_i(t)}{\sigma_{X}(t)}\left(1-(1 \pm \varepsilon)\frac{\sigma_{X}(T)-\sigma_{X}(t)}{\sigma_{X}(T)}\right),\ \ t\ge0,\ 1\le i\le n.$$
Since the analysis of $\pi_{+\vn}(u)$ and $\pi_{-\vn}(u)$ are the same, next we only discuss $\pi_{+\vn}(u)$ for fixed $\vn\in(0,1)$.
The variance function $\sigma_Y(t)$ of $Y_{+\vn,1 }(t)$ attains its maximum over $[\delta,T]$ at unique point $T$ with
\BQNY
\sigma_Y(t)=1-A(1+\vn)(T-t)^\mu(1+\oo(1)), \ \ \text{as}\ t\rw T.
\EQNY
Further, by \eqref{eq:rX}
\BQN
r_Y(s,t)=Corr(Y_{+\vn,1 }(s), Y_{+\vn,1 }(t))=1- D|t-s|^{\nu}+o(|t-s|^{\nu}), \quad \min(t,s)\rw T.
\EQN
Moreover, in view of Assumption A2 for $s,t\in[\delta,T]$ 
\begin{eqnarray*}
\E{(Y_{+\ve,1}(t)-Y_{+\ve,1}(s))^{2}}
&\leq & \mbbC\ |s-t|^{\gamma}.
\end{eqnarray*}
Consequently, by \netheo{MainThm2}
\BQNY
\pi_{+\vn}(u) 
=  \WW_{\nu,\mu}^\vn \LT(m_u(T)\RT)^{(2/\nu-2/\mu)_+} \Upsilon_n \LT(m_u(T)\RT)
\ooo, \quad u\to \IF,
\EQNY
where
$$
\WW_{\nu,\mu}^\vn=\left\{
            \begin{array}{ll}
D^{1/\nu} \frac{\Gamma(1/\mu+1)}{((1+\vn)A)^{1/\mu}}\mathcal{H}_{\nu}, & \hbox{if } \nu<\mu ,\\
\P_{\nu,\mu}^{ A(1+\vn) D^{-1}}, & \hbox{if } \nu=\mu,\\
1&  \hbox{if } \nu>\mu.
              \end{array}
            \right.
$$
Letting $\vn\rw0$, we conclude that
\BQNY
\Pi(u)= \WW_{\nu,\mu}^0 
\LT(m_u(T)\RT)^{(2/\nu-2/\mu)_+}
\Upsilon_n \LT(m_u(T)\RT)\ooo
\EQNY
as $u\to \infty$.
In addition, by the Borell-TIS inequality, for $u$ sufficiently large
\BQNY
\pk{\sup_{t\in[0,\delta]} \xih >u}&\le& \pk{\sup_{t\in[0,\delta]}\chi_n(t)>u}\nonumber\\
&=&\pk{\sup_{(t,\vk{v})\in \mathcal{G}_{\delta}}Y(t,\vk{v})>u}\nonumber\\
&\le&\exp\LT(-\frac{(u-\mathbb{C})^2}{2\sigma_X^2(\delta)}\RT),
\EQNY
and thus the claim follows from the last two formulas. \QED

\def\BkS{B(k,S_1,u)}
\section{Appendix}
This section is dedicated to the proof of \eqref{eq:doubleS1}. Let
$$
\BkS=u+c(kS_1u^{-2/\alpha})^\beta,\ \  k\inn,\ S_1>0,\ u>0.
$$
The double sum $\Sigma_{\chi}(u)$ can be divided into two parts, i.e.,
$$\Sigma_{\chi}(u)=\sum_{0\le k<j\le h(u)-1} \pk{\sup_{t\in\det_k} \chicb >u, \sup_{t\in\det_j} \chicb >u }=:\Sigma_{\chi,1}(u)+\Sigma_{\chi,2}(u),$$
 where  $\Sigma_{\chi,1}(u)$ is the sum \re{for}  $j=k+1$, and $\Sigma_{\chi,2}(u)$ is the sum \re{for}  $j>k+1$. We first give the estimation of the first sum.
It follows that
\BQN\label{eq:Sig1}
\Sigma_{\chi,1}(u)
&\le&\sum_{ k=0}^{ h(u)} \pk{\sup_{t\in\det_k}\chi_n(t)>\BkS , \sup_{t\in\det_{k+1}}\chi_n(t)>\BkS }.
\EQN
Further, we have that
\BQNY
 &&\pk{\sup_{t\in\det_k}\chi_n(t)>\BkS , \sup_{t\in\det_{k+1}}\chi_n(t)>\BkS }\\
&=& \pk{\sup_{t\in\det_k}\chi_n(t)>\BkS}+\pk{\sup_{t\in\det_{k+1}}\chi_n(t)>\BkS }\\
&&-\pk{\sup_{t\in\det_k\cup\det_{k+1}}\chi_n(t)>\BkS },
\EQNY
which, in the light of the reasoning of  \eqref{eq:upper1} gives that
$$
\lim_{S_1\rw\IF}\limsup_{u\rw\IF}\frac{\Sigma_{\chi,1}(u)}{A_2(u)}\le \mbbC \lim_{S_1\rw\IF}\frac{2\Ha[0,S_1]-\Ha[0,2S_1]}{S_1}=0.
$$
In order to estimate $\Sigma_{\chi,2}(u)$, we introduce
 a Gaussian \re{random} field
$$
Y(t,\vk{v})=\sum_{i=1}^n v_i X_i(t), \quad t\ge0,\ \vk{v}=(v_1,\cdots,v_n)\in\R^n.
$$
In the light of Piterbarg (1996) 
$$
\underset{t\in[0,S_1]}\sup \chi_n(t)=\underset{(t,\vk{v})\in\mathcal{G}_{S_1}}\sup Y(t,\vk{v}),
$$
where  $\mathcal{G}_{S_1}=[0,S_1]\times\mathcal{S}_{n-1}$, with $\mathcal{S}_{n-1}$ being the unit sphere in $\R^n$. Consequently,
\BQN\label{eq:Sig1}
\Sigma_{\chi,2}(u)
\le\sum_{k=0}^{h(u)-1}\sum_{j=k+2}^{h(u)-1} \pk{\underset{(t,\vk{v})\in \det_k\times\mathcal{S}_{n-1}}\sup Y(t,\vk{v})>\BkS , \underset{(t,\vk{v})\in \det_{j}\times\mathcal{S}_{n-1}}\sup Y(t,\vk{v})>\BkS }.
\EQN
We split the sphere $\On$ into sets of small diameters $\{\parO_i, 0\le i\le \NN\}$,  
where
$$
\NN=\sharp\{\parO_i\}<\IF.
$$
Further, we see that the summand on the right-hand side of \eqref{eq:Sig1} is not greater than $\Sigma_1^{k,j}(u)+\Sigma_2^{k,j}(u)$,
with
\BQNY
\Sigma_1^{k,j}(u)&=& \underset{\parO_i\cap\parO_l=\emptyset}{\sum_{0\le i,l\le \NN}}\pk{\underset{(t,\vk{v})\in \det_k\times\parO_i}\sup Y(t,\vk{v})>\BkS , \underset{(t,\vk{v})\in \det_j\times\parO_l}\sup Y(t,\vk{v})>\BkS }\\
\Sigma_2^{k,j}(u)&=& \underset{\parO_i\cap\parO_l\neq\emptyset}{\sum_{0\le i,l\le \NN}}\pk{\underset{(t,\vk{v})\in \det_k\times\parO_i}\sup Y(t,\vk{v})>\BkS , \underset{(t,\vk{v})\in \det_j\times\parO_l}\sup Y(t,\vk{v})>\BkS },
\EQNY
where  $\parO_i\cap\parO_l\neq\emptyset$ means $\parO_i, \parO_l$ are identical or adjacent, and $\parO_i\cap\parO_l=\emptyset$ means $\parO_i, \parO_l$ are neither identical nor adjacent. Denote the distance of two sets $\vk{A},\vk{B}\in\R^n,  n\inn,$ as
$$
\rho(\vk{A},\vk{B})=\underset{\vk{x}\in \vk{A}, \vk{y}\in \vk{B}}\inf||\vk{x}-\vk{y}||^2.
$$
If $\parO_i\cap\parO_l=\emptyset$ then there exists some small positive constant $\rho_0$ (independent of $i,l$) such that $\rho(\parO_i,\parO_l)>\rho_0$. Next, we  estimate $\Sigma_1^{k,j}(u)$.
\COM{
Let
$$
Z(t,\vk{v},s,\vk{w})=Y_{1}(t,\vk{v})+Y_{2}(s,\vk{w}),\  t,s\in [0,\IF),\ \vk{v}, \vk{w}\in\R^n,
$$
with
$$
Y_{1}(t,\vk{v}):=\sum_{i=1}^n v_i X_i(t), \ \
Y_{2}(s,\vk{w}):=\sum_{i=1}^n w_i X_i(s), \quad t,s\in[0,\IF),\ \vk{v},\vk{w}\in\On.
$$
}
For any $u\ge 0$
\BQNY
&&\pk{\underset{(t,\vk{v})\in \det_k\times\parO_i}\sup Y(t,\vk{v})>\BkS , \underset{(t,\vk{v})\in \det_j\times\parO_l}\sup Y(t,\vk{v})>\BkS }\\
&\le&\pk{\underset{\vk{v}\in\parO_i,\vk{w}\in\parO_l}{\sup_{(t,s)\in\det_k\times\det_{j}}}Z(t,\vk{v},s,\vk{w})>2u },
\EQNY
where
$$
Z(t,\vk{v},s,\vk{w})=Y(t,\vk{v})+Y(s,\vk{w}),\  t,s\ge0,\ \vk{v}, \vk{w}\in\R^n.
$$
When $u$ is sufficiently large for $ (t,s)\in\det_k\times\det_j\subset[0,1]^2, \vk{v}\in\parO_i\subset[-2,2]^n,\vk{w}\in\parO_l\subset[-2,2]^n,$ with $\rho(\parO_i,\parO_l)> \rho_0$ we have
\BQNY
Var(Z(t,\vk{v},s,\vk{w}))
&=&4-\LT(2(1-r(s-t))+r(s-t)||\vk{v}-\vk{w}||^2\RT)\\
&\le&4(1-\delta_0),
\EQNY
for some $\delta_0>0$.  Therefore, it follows from Borell-TIS inequality (see e.g., Adler and Taylor (2007)) that
\BQNY
\Sigma_1^{k,j}(u)\le \mbbC \NN  \exp\LT(-\frac{(u-a)^2}{2(1-\delta_0)}\RT),\  \text{with} \  a=\E{\underset{(\vk{v},\vk{w})\in[-2,2]^{2n}}{\underset{(t,s)\in[0,1]^2}\sup}Z(t,\vk{v},s,\vk{w})}<\IF.
\EQNY
Consequently,
\BQN\label{eq:Sig1_2}
\limsup_{u\rw\IF}\frac{\sum_{k=0}^{h(u)-1}\sum_{j=k+2}^{h(u)-1} \Sigma_1^{k,j}(u)}{A_2(u)}  =0 .
\EQN

Next, we estimate $\sum_{k=0}^{h(u)-1}\sum_{j=k+2}^{\IF}\Sigma_2^{k,j}(u)$.
The stationarity of  $\{Y(t,\vk{v}),t\ge0, \vk{v}\in\On\}$ implies
\BQNY
\sum_{k=0}^{h(u)-1}\sum_{j\ge k+2}\Sigma_2^{k,j}(u)&\le& \mbbC\sum_{k=0}^{h(u)-1}\sum_{j\ge2} \pk{\underset{(t,\vk{v})\in \det_0\times\parO_i}\sup Y(t,\vk{v})>\BkS , \underset{(t,\vk{v})\in \det_j\times\parO_l}\sup Y(t,\vk{v})>\BkS }
\EQNY
for some fixed $\parO_i,\parO_l$ satisfying $\parO_i\cap\parO_l\neq\emptyset$. Additionally, diam$(\parO_i\cup\parO_l)$ can be chosen sufficiently small such that $\parO_i,\parO_l$ are in $\pO_0$, which is a subset of $\On$ and includes $(1,0,\cdots,0)$, and further on $\pO_0$ we can find a one-to-one projection $g$ from it to the corresponding points where the first component is 1, i.e., $g\vk{v}=(1,v_2,\cdots,v_n)$ for all $\vk{v}=(v_1,v_2,\cdots,v_n)\in\pO_0$.

 Let
$$ \wtre_0=\LT[0,\frac{S_2}{u}\RT]^{n-1}\cap g\pO_0,\  \wtre_{\vk{k}}=\prod_{i=1}^{n-1}\LT[k_i\frac{S_2}{u},(k_i+1)\frac{S_2}{u}\RT]\cap g\pO_0,\  \vk{k}\in \mathbb{Z}^{n-1}$$
\re{and}
$$\CK_i=\{\vk{k}:\wtre_{\vk{k}}\cap g\parO_i\neq \emptyset\},\ \ \ \CK_l=\{\vk{k}:\wtre_{\vk{k}}\cap g\parO_l\neq \emptyset\}.$$ 
With these notation, we have that
\BQNY
\sum_{k=0}^{h(u)-1}\sum_{j\ge k+2}\Sigma_2^{k,j}(u)&\le& \mbbC\sum_{k=0}^{h(u)-1}\sum_{j\ge2}\sum_{\vk{i}\in\CK_i}\sum_{\vk{l}\in\CK_l}\\
 &&\pk{\underset{(t,\vkt{v})\in \det_0\times\wtre_{\vk{i}}}\sup Y(t,\vkt{v})>\BkS , \underset{(t,\vkt{v})\in \det_j\times\wtre_{\vk{l}}}\sup Y(t,\vkt{v})>\BkS }.
\EQNY
The last sums on the right-hand side can be divided into two terms $I_i(u), i=1,2,$ according to whether $\wtre_{\vk{i}}\cap\wtre_{\vk{l}}\neq\emptyset$ or not. We derive that
\def\vki{\vk{i}}
\def\vkl{\vk{l}}
\BQNY
&&\pk{\underset{(t,\vkt{v})\in \det_0\times\wtre_{\vki}}\sup Y(t,\vkt{v})>\BkS , \underset{(t,\vkt{v})\in \det_j\times\wtre_{\vkl}}\sup Y(t,\vkt{v})>\BkS }\\
&\le&\pk{\underset{\vkt{v}\in\wtre_{\vki},\vkt{w}\in\wtre_{\vkl}}{\sup_{(t,s)\in\det_0\times\det_{j}}}Z(t,\vkt{v},s,\vkt{w})>2\BkS },
\EQNY
where
$$
Z(t,\vkt{v},s,\vkt{w})=Y(t,\vkt{v})+Y(s,\vkt{w}),\  t,s\ge0,\ \vkt{v}, \vkt{w}\in\R^{n-1}.
$$

\COM{
Notice that the summand in $\Sigma_2^{k,j}$ can be rewritten as
\BQNY
&&\pk{\underset{(t,\vk{s})\in \det_0\times\parO_i}\sup Y(t,\vk{s})>u+c(kSu^{-2/\alpha})^\beta , \underset{(t,\vk{s})\in \det_{j-k}\times\parO_l}\sup Y(t,\vk{s})>u+c(jSu^{-2/\alpha})^\beta }\\
&\le&\pk{\underset{\vk{v}\in\parO_i}{\sup_{(t,s)\in\det_0\times\det_{j-k}}}Z(t,\vk{v},s,\vk{v})>2u+2c(kSu^{-2/\alpha})^\beta, \underset{\vk{v}\in\parO_l}{\sup_{(t,s)\in\det_0\times\det_{j-k}}}Z(t,\vk{v},s,\vk{v})>2u+2c(kSu^{-2/\alpha})^\beta},
\EQNY
where in the last inequality we used the facts that the trajectories of $\{Z(t,\vk{v},s,\vk{v})\}$ is continuous in $\vk{v}$ and that
 $\rho(\parO_i,\parO_l)\le\rho_0$ (here $\rho_0$ is sufficiently small).
Notice that $\{Z(t,\vk{v},s,\vk{v})\}$ is stationary in $\vk{v}$ on $\On$. We can, without loss of generality, assume that $\parO_i,\parO_l$ is on the sphere near the point $(1,0,\cdots,0)$.\\
 Consequently,
\BQN\label{eq:Sig1_3}
\Sigma_2^{k,j}&\le& \underset{\rho(\parO_i,\parO_l)\le \rho_0}{\sum_{\parO_i,\parO_l}}(P_u^{k,j}(\parO_i)+P_u^{k,j}(\parO_l)-P_u^{k,j}(\parO_i\cup\parO_l)),
\EQN
where
\BQNY
P_u^{k,j}(D)&=&\pk{\underset{\vk{v}\in D}{\sup_{(t,s)\in\det_0\times\det_{j-k}}}Z(t,\vk{v},s,\vk{v})>2u+2c(kSu^{-2/\alpha})^\beta}.
\EQNY
Hence it follows from \eqref{eq:Sig1_2} and \eqref{eq:Sig1_3} that
\BQNY
\Sigma_{\chi,u}
&\le&\sum_{0\le k<j\le h(u)-1}K \NN u^{-1}\exp\LT(-\frac{(u-a)^2}{2(1-\delta_0)}\RT)\\
&&+\underset{\rho(\parO_i,\parO_l)\le \rho_0}{\sum_{\parO_i,\parO_l}}\sum_{0\le k<j\le h(u)-1}(P_u^{k,j}(\parO_i)+P_u^{k,j}(\parO_l)-P_u^{k,j}(\parO_i\cup\parO_l)),
\EQNY
from which we see that
\BQNY
\Sigma_{\chi,u}=\oo\LT(
u^{2/\alpha-1/\beta+n-2}
\exp\left(-\frac{u^2}{2}\right)\RT), \ \ \ u\rw\IF,
\EQNY
once we obtain
\BQN\label{eq:P_u}
\sum_{0\le k<j\le h(u)-1}P_u^{k,j}(\parO_i)=\oo\LT(
u^{2/\alpha-1/\beta+n-2}
\exp\left(-\frac{u^2}{2}\right)\RT), \ \ \ u\rw\IF,
\EQN
for any $\parO_i$ with $||\parO_i||\le\delta$ on the sphere  near the point $(1,0,\cdots,0)$. The right-hand side of Eq. \eqref{eq:P_u} can be rewritten as
\BQN\label{eq:P_u_1}
\sum_{0\le k<j\le h(u)-1}P_u^{k,j}(\parO_i)&=& \sum_{k=0} ^{h(u)-1}\sum_{j=k+1}P_u^{k,j}(\parO_i)+\sum_{k=0} ^{h(u)-1}\sum_{j>k+1}P_u^{k,j}(\parO_i)
\EQN
}


It follows that, for $(t,s)\in\det_0\times\det_{j}$, $\vkt{v}\in\wtre_{\vki},\vkt{w}\in\wtre_{\vkl}$, diam$(\pO_0)$ sufficiently small, and $u$ sufficiently large
\BQN\label{eq:VZ}
2\le Var(Z(t,\vkt{v},s,\vkt{w}))
\le 4\LT(1-\frac{1}{4} ((j-1)S_1)^\alpha u^{-2}\RT).
\EQN
Further, set $\overline{Z}(t,\vkt{v},s,\vkt{w})=Z(t,\vkt{v},s,\vkt{w})/\sqrt{Var(Z(t,\vkt{v},s,\vkt{w}))}$. Borrowing the arguments of the proof of Lemma 6.3 in Piterbarg (1996) we show that
\BQNY
\E{\overline{Z}(t,\vkt{v},s,\vkt{w})-\overline{Z}(t',\vkt{v}',s',\vkt{w}')}^2\le 4\LT(\E{Y(t,\vkt{v})-Y(t',\vkt{v}')}^2+\E{Y(s,\vkt{w})-Y(s',\vkt{w}')}^2\RT).
\EQNY
Moreover, as in  Lemma 10 of Piterbarg (1994b),  for diam$(\pO_0)$ sufficiently small, and $u$ sufficiently large,
\BQNY
\E{Y(t,\vkt{v})-Y(t',\vkt{v}')}^2
&\le&4  |t-t'|^\alpha+2\sum_{i=2}^n(v_i-v'_i)^2.
\EQNY
Therefore 
\BQN\label{eq:Zbar}
\E{\overline{Z}(t,\vkt{v},s,\vkt{w})-\overline{Z}(t',\vkt{v}',s',\vkt{w}')}^2&\le& 16  |t-t'|^\alpha+16 |s-s'|^\alpha+8\sum_{i=2}^n(v_i-v'_i)^2+8\sum_{i=2}^n(w_i-w'_i)^2\nonumber\\
&\le&2(1-r_\zeta(|t-t'|,|s-s'|,\vkt{v}-\vkt{v}',\vkt{w}-\vkt{w}')),
\EQN
where
$$r_\zeta(t,s,\vkt{v},\vkt{w})=\exp\LT(-9  t^\alpha-9  s^\alpha-5\sum_{i=2}^nv_i^2-5\sum_{i=2}^nw_i^2\RT), \ \  t,s\ge0, \vkt{v},\vkt{w}\in\R^{n-1}
$$
is the covariance function of a stationary Gaussian \re{random} field $\{\zeta(t,s,\vkt{v},\vkt{w}), t,s\ge0, \vkt{v},\vkt{w}\in\R^{n-1}\}$. Consequently, in view of \eqref{eq:VZ} and \eqref{eq:Zbar}, and thanks to Slepian's Lemma, we obtain
\BQNY
&&\pk{\underset{(t,\vkt{v})\in \det_0\times\wtre_{\vki}}\sup Y(t,\vkt{v})>\BkS , \underset{(t,\vkt{v})\in \det_j\times\wtre_{\vkl}}\sup Y(t,\vkt{v})>\BkS }\\
&\le&\pk{\underset{\vkt{v}\in\wtre_{\vki},\vkt{w}\in\wtre_{\vkl}}{\sup_{(t,s)\in\det_0\times\det_{j}}}\zeta(t,s,\vkt{v},\vkt{w})>
\frac{2\BkS}{\sqrt{4- ((j-1) S_1)^\alpha u^{-2}}} }.
\EQNY
Since, for any cube $\wtre_{\vki}$ in $\R^{n-1}$ there are $3^{n-1}$ cubes having non-empty intersection with it,  we have
\BQN\label{eq:I1u1}
I_1(u)
&\le&\sum_{k=0}^{h(u)-1}\sum_{j\ge2}\underset{\wtre_{\vk{i}}\cap\wtre_{\vk{l}}\neq\emptyset}{\sum_{\vk{i}\in\CK_i}\sum_{\vk{l}\in\CK_l}}\pk{\underset{\vkt{v}\in\wtre_{\vki},\vkt{w}\in\wtre_{\vkl}}{\sup_{(t,s)\in\det_0\times\det_{j}}}\zeta(t,s,\vkt{v},\vkt{w})>
\frac{2\BkS}{\sqrt{4- ((j-1) S_1)^\alpha u^{-2}}} }\nonumber\\
&\le&3^{n-1}\sum_{k=0}^{ h(u)-1}\sum_{j\ge2}\sum_{\vk{i}\in\CK_i}\pk{\underset{\vkt{v}\in\wtre_{\vki},\vkt{w}\in\wtre_{\vkl}}{\sup_{(t,s)\in\det_0\times\det_{j}}}
\zeta(t,s,\vkt{v},\vkt{w})>
\frac{2\BkS}{\sqrt{4- ((j-1) S_1)^\alpha u^{-2}}} },
\EQN
with some $\wtre_{\vkl}$ adjacent or identical with $\wtre_{\vki}$.
It follows further from \netheo{ThmlamS} that
\BQNY
\pk{\underset{\vkt{v}\in\wtre_{\vki},\vkt{w}\in\wtre_{\vkl}}{\sup_{(t,s)\in\det_0\times\det_{j}}}\zeta(t,s,\vkt{v},\vkt{w})>
\frac{2\BkS}{\sqrt{4- ((j-1) S_1)^\alpha u^{-2}}} }&\le&\LT(\Ha[0,9^{\frac{1}{\alpha}}S_1]\RT)^{2} \LT(\mathcal{H}_2[0,\sqrt{5}S_2]\RT)^{2(n-1)}\\
&& \frac{1}{\sqrt{2\pi}u}\exp\LT(-\frac{4\BkS^2}{2\LT(4- ((j-1) S_1)^\alpha u^{-2}\RT)}\RT)\ooo
\EQNY
as $u\rw\IF.$ Inserting the last formula into \eqref{eq:I1u1} and noting that
$$
\sharp\{\CK_i\}= V_{n-1}(g\parO_i)S_2^{-(n-1)}u^{n-1}\ooo,\ \  \text{as}\ u\rw\IF
$$
we derive that
\BQNY
I_1(u)&\le&\mbbC S_1^2S_2^{n-1}\sum_{k=0}^{h(u)-1}\sum_{j\ge2}
\frac{1}{\sqrt{2\pi}}u^{n-2}\exp\LT(-\frac{u^2}{2}-c(kS_1u^{\frac{1}{\beta}-\frac{2}{\alpha}})^\beta-\frac{1}{8}((j-1) S_1)^\alpha \RT).
\EQNY
Thus, in the light of the reasoning of \eqref{eq:upper1}, we conclude that
\BQN\label{eq:I1u}
\limsup_{S_1\rw\IF}\limsup_{u\rw\IF}\frac{I_1(u)}{A_2(u)}\le \mbbC \limsup_{S_1\rw\IF}S_1 S_2^{n-1}\exp\LT(-\frac{1}{8}S_1^\alpha\RT)=0.
\EQN
Moreover, in view of the reasoning of \eqref{eq:VZ}, when $\wtre_{\vk{i}}\cap\wtre_{\vk{l}}=\emptyset$, we obtain
\BQNY
2\le Var(Z(t,\vkt{v},s,\vkt{w}))
\le 4- ((j-1) S_1)^\alpha u^{-2}-||\vk{l}-\vk{i}||^2S_2^2u^{-2}
\EQNY
and thus
\BQNY
I_2(u)
&\le&\sum_{k=0}^{h(u)-1}\sum_{j\ge2}\underset{\wtre_{\vk{i}}\cap\wtre_{\vk{l}}=\emptyset}{\sum_{\vk{i}\in\CK_i}\sum_{\vk{l}\in\CK_l}}\pk{\underset{\vkt{v}\in\wtre_{\vki},\vkt{w}\in\wtre_{\vkl}}{\sup_{(t,s)\in\det_0\times\det_{j}}}\zeta(t,s,\vkt{v},\vkt{w})>
\frac{2\BkS}{\sqrt{4- ((j-1) S_1)^\alpha u^{-2}-||\vk{i}-\vk{l}||^2S_2^2u^{-2}}} }\\
&\le& \sum_{k=0}^{h(u)-1}\sum_{j\ge2}\sum_{\vk{i}\in\CK_i}\underset{\vk{l}\neq0}{\sum_{\vk{l}\in\R^{n-1}}}
\pk{\underset{\vkt{v}\in\wtre_{\vk{0}},\vkt{w}\in\wtre_{\vkl}}{\sup_{(t,s)\in\det_0\times\det_{j}}}\zeta(t,s,\vkt{v},\vkt{w})>
\frac{2\BkS}{\sqrt{4- ((j-1) S_1)^\alpha u^{-2}-||\vk{l}||^2S_2^2u^{-2}}} }.
\EQNY
Similar to \eqref{eq:I1u}, we conclude that
\BQNY
\limsup_{S_1\rw\IF}\limsup_{u\rw\IF}\frac{I_2(u)}{A_2(u)} =0,
\EQNY
hence \eqref{eq:doubleS1} follows.\\

\textbf{Acknowledgments.} Partial support from the Swiss National Science Foundation Project 200021-1401633/1 and by the project RARE -318984
 (a Marie Curie International Research Staff Exchange Scheme Fellowship within the 7th European Community Framework Programme) is kindly acknowledged.


\begin{thebibliography}{100} \small
\bibitem{}Adler, R.J. and Taylor, J.E., 2007. {\it Random Fields and Geometry}. Springer.
\bibitem{} Albin, J.M.P., On  extremal theory for stationary processes. Ann. Probab. 18 (1990), 92-128.
\bibitem{}
Albin, J.M.P. and Jaru$\check{s}$kov\'{a}, D., On a test statistic for linear trend. Extremes 6 (2003), 247-258.
\bibitem{}Aronowich, M. and Adler, R. J., Behaviour of $\chi^2$ processes at extrema.
Advances in Applied Probability  17 (1985),  280-297.

\bibitem{} Belyaev, Yu. K. and Nosko, V.P., Characteristics of exursions above a high level for a Gaussian process and its envelope. Theory Probab. Appl. 13 (1969), 298-302.
\bibitem{} Berman, M.S. 1992. {\it Sojourns and Extremes of Stochastic Processes}. Wadsworth and Brooks/ Cole, Boston.


\bibitem{}
Bojdecki, T., Gorostiza, L., and Talarczyk, A.,  Sub-fractional Brownian motion and its relation to
occupation times, Statistics and Probability Letters 69 (2004),  405-419.


\bibitem{17} D\c{e}bicki, K.,  Ruin probability for Gaussian integrated processes. Stoch. Proc. Appl.
     98 (2002), 151-174.

\bibitem{} D\c{e}bicki, K., Sikora, G., Finite time asymptotics of fluid and ruin models:
multiplexed fractional Brownian motions case. Applicationes Mathematicae 38 (2011),  107-116.

\bibitem{} D\c{e}bicki, K. and Tabi\'{s}, K., Extremes of time-average stationary Gaussian processes. Stoch. Proc. Appl. 121 (2011), 2049-2063.
\bibitem{} Dieker, A.B. and Yakir, B.,  On asymptotic constants in the theory of Gaussian processes. Bernoulli. to appear. (2013).

\bibitem{falk}
Falk, M., H\"usler, J.,  and Reiss, R.D.,  2010. {\it Laws of Small
Numbers: Extremes and Rare Events.} DMV Seminar Vol. {23}, 2nd edn., Birkh\"auser, Basel.

%


\bibitem{}
Hashorva, E., Kabluchko, Z., and W\"ubker, A.,  Extremes of independent chi-square random vectors.
Extremes 15 (2012), 35--42.

\bibitem{} Houdr\'{e}, C. and Villa, J., An example of infinite dimensional quasi-helix\cH{.} Contemporary Mathematics, American Mathematical Society 336 (2003), 195-201.

\bibitem{} Jaru\v{s}kov\'{a}, D., Asymptotic behaviour of a test statistic for detection of change in
mean of vectors. J. Stat. Plan. Inf. 140 (2010), 616--625.

\bibitem{}
Jaru$\check{s}$kov\'{a}, D. and  Piterbarg, V.I.,  Log-likelihood ratio test for detecting transient
change. Stat. Probab. Lett. 81 (2011), 552-559.

\bibitem{}
Kabluchko, Z.,  Extremes of independent Gaussian processes.  Extremes 14  (2011), 285--310.


\bibitem{} Kozachenko, Y. and Moklyachuk, O., Large deviation probabilities for square-Gaussian stochastic processes. Extremes 2 (1999), 269-293.

\bibitem{} Lindgren, G., Extreme values and crossings for the $\chi^{2}$-process and other functions of multidimensional Gaussian proceses with reliability applications. Adv. Appl. Probab. 12 (1980a), 746-774.
\bibitem{} Lindgren, G., Point processes of exits by bivariate Gaussian processes and extremal theory for the $\chi^{2}$-process and its concominats. J. Multivar. Anal. 10 (1980b), 181-206.

\bibitem{} Lindgren, G., Slepian models for $\chi^{2}$-process with dependent components with application to envelope upcrossings. J. Appl. Probab. 26 (1989), 36-49.
    \bibitem{} Michna, Z., Remarks on Pickands theorem. Preprint. (2009).  Available at http://arxiv.org/abs/0904.3832

\bibitem{} Pickands III, J., Upcrossing probabilities for stationary Gaussian processes. Trans. Amer. Math. Soc.
145 (1969a), 51-73.
\bibitem{}
Pickands III, J., Asymptotic properties of the maximum in a stationary Gaussian process. Trans. Amer. Math. Soc. 145 (1969b), 75-86.
\bibitem{}
Piterbarg, V.I., On the paper by J. Pickands "Upcrosssing probabilities for stationary Gaussian processes", Vestnik
Moscow. Univ. Ser. I Mat. Mekh. 27 (1972), 25-30. English transl. in Moscow Univ. Math. Bull., 27 (1972).

\bibitem{} Piterbarg, V.I., High deviations for multidimensional stationary Gaussian processes with independent components. In: Zolotarev,V.M.(Ed.), Stability Problems for Stochastic Models (1994a), 197-210.

\bibitem{} Piterbarg, V.I., High exsursions for nonstationary generalized chi-square processes. Stoch. Proc. Appl. 53 (1994b), 307-337.


\bibitem{} Piterbarg, V.I., 1996. {\it Asymptotic Methods in the Theory of Gaussian Processes and Fields}. In: Transl. Math. Monographs, vol. 148. AMS, Providence, RI.

\bibitem{} Piterbarg, V.I.,  Large deviations of a storage process with fractional Browanian motion as input. Extremes 4 (2001), 147-164.

\bibitem{} Piterbarg, V.I., Prisyazhnyuk, V., Asymptotic behavior of the probability of a large excursion for a
nonstationary Gaussian processes. Teor. Veroyatnost. i Mat. Statist. 18 (1978), 121-133.

\end{thebibliography}
\end{document}